\input amstex
\documentstyle{amsppt}
\vsize=7.4in 
\hsize=5.5in 
\nologo

\def\d{\Delta}
\def\p{\partial}
\def\n{\nabla}
\def\la{\langle}
\def\ra{\rangle}

\def\-{\setminus}

\topmatter
\title
Connectedness at infinity of complete K\"ahler manifolds and 
locally symmetric spaces
\endtitle

\rightheadtext{Connectedness at infinity}
\leftheadtext{Peter Li and Jiaping Wang}
\author
Peter Li and Jiaping Wang
\endauthor
\affil
University of California, Irvine\\
University of Minnesota
\endaffil

\footnotetext" "{The first author was partially supported 
by NSF grant DMS-0503735.
The second author was partially supported by NSF grant DMS-0404817.}

\address
Department of Mathematics, University of California,
Irvine, CA 92697-3875
\endaddress

\email
pli\@math.uci.edu
\endemail

\address
School of Mathematics, University of Minnesota,
Minneapolis, MN 55455
\endaddress

\email
jiaping\@math.umn.edu
\endemail

\date
February 13, 2007
\enddate

\abstract
One of the main purposes of this paper is to prove that on a complete K\"ahler 
manifold of dimension $m$, if the holomorphic bisectional curvature is bounded 
from below by -1 and the minimum spectrum $\lambda_1(M) \ge m^2$,  then it must 
either be connected at infinity or isometric to $\Bbb R \times N$ with a specialized metric, 
with $N$ being compact.  Similar type results are also
proven for irreducible, locally symmetric spaces of noncompact type.  Generalizations
to complete K\"ahler manifolds satisfying a weighted Poincar\'e inequality are also considered.
\endabstract 
\endtopmatter

\document
\subheading
{Table of Contents}

\noindent
\S0 Introduction \newline
\noindent
\S1 K\"ahler manifolds with maximum $\lambda_1$\newline
\noindent
\S2 Locally Symmetric Spaces\newline
\noindent
\S3 Nonparabolic Ends \newline
\noindent
\S4 Parabolic Ends with $\rho \to 0$

\heading
\S0 Introduction
\endheading

One of the purposes of the current paper is to continue our study of
K\"ahler manifolds whose holomorphic bisectional curvature is 
bounded from below by $-1$, while the greatest lower bound of its 
spectrum is bounded from below by a positive constant. For a  
complete K\"ahler manifold $M^m$ of complex dimension $m,$  
the holomorphic bisectional curvature $BK_M$  is said to have 
the lower bound $BK_M \ge -1$ if

$$
R_{i\bar i j \bar j} \ge -(1 + \delta_{ij})
$$
for all unitary frames $\{e_1, e_2 \dots, e_m\}.$  Under the 
assumption that $BK_M \ge -1$, the authors proved in \cite{L-W3}
a comparison theorem for the Laplacian of the distance function, 
and as a consequence (Corollary 1.10 of \cite{L-W3}), they 
established that

$$
\lambda_1(M) \le m^2,
$$
where $\lambda_1(M)$ is the greatest lower bound for the 
$L^2$-spectrum of the Laplacian acting on functions.  
This theorem can be viewed as an  analogue to Cheng's 
theorem \cite{C} for K\"ahler manifolds.  In the same paper, 
they also considered the equality case of this estimate 
on $\lambda_1(M)$ and obtained some partial results concerning
the connectedness at infinity for complete
K\"ahler surfaces. Here, we prove a more complete result that holds for all dimensions.

\proclaim{Theorem A}  Let $M^m$ be a  complete K\"ahler manifold 
of complex dimension $m\ge 2$ with holomophic bisectional curvature 
bounded by

$$
BK_M\ge -1.
$$ 
If  $\lambda_1(M) = m^2,$ then either
\roster
\item $M$ is connected at infinity; or
\item $M$ is diffeomorphic to $\Bbb R \times N$ with
$N$ being a compact compact.  Moreover, 
$M$ is covered by the complex hyperbolic space $\Bbb {CH}^m$
if $M$ has bounded curvature.
\endroster
\endproclaim

In the process, we also show the following theorem.

\proclaim{Theorem B} Let $M^m$ be a complete K\"ahler manifold with
its Ricci curvature satisfying

$$
\text{Ric}_M \ge -2(m+1).
$$ 
If $\lambda_1(M)>\frac{m+1}{2},$
then $M$ must have only one infinite volume end.
\endproclaim

We would like to point out that according to our normalization 
in \cite{L-W3}, the assumption that $BK_M \ge -1$ implies 
that $\text{Ric}_M \ge -2(m+1).$ Also,
the lower bound $\frac{m+1}2$ for $\lambda_1(M)$ in the theorem
is best possible as
one can find counterexample of the form $M = \Sigma \times N,$ 
where $N^{m-1}$ is a compact K\"ahler manifold and $\Sigma$ 
a complete Riemann surface with more than one infinite volume ends.  
These examples motivate us to propose the following conjecture.

\proclaim{Conjecture C} Let $M^m$ be a complete K\"ahler manifold 
with Ricci curvature satisfying
$$
\text{Ric}_M \ge -2(m+1)
$$
and
$$
\lambda_1(M) \ge \frac{m+1}2.
$$
Then either
\roster
\item $M$ has only one infinite volume end; or
\item it is the total space of a holomorphic fiber bundle 
$N \rightarrow M \rightarrow \Sigma$, with compact totally geodesic 
fiber $N$, over a complete Riemann surface $\Sigma.$  Moreover, 
$\lambda_1(M) = \frac{m+1}2.$
\endroster
\endproclaim 

While these results were motivated by the authors earlier studies, 
\cite{L-W1} and \cite{L-W2}, on Riemannian manifolds where they 
gave a rather complete picture of similar type theorems, it should be
noted that different type of arguments are required here. 
It turns out the new approaches can be adapted to deal with 
more general manifolds.
In particular, we have an analogous result to Theorem A for locally
symmetric spaces.
 
\proclaim{Theorem D} Let $M$ be an irreducible locally symmetric space 
covered by an irreducible symmetric space of noncompact type 
$\tilde M = G/K.$  Suppose $\lambda_1(M) = \lambda_1(\tilde M).$  
Then either 
\roster
\item $M$ is connected at infinity; or
\item it is isometric to $\Bbb R \times N$ with $N$ being
a compact quotient of the horosphere of $\tilde M$.  
\endroster
\endproclaim

Many of the aforementioned results can be generalized to complete 
K\"ahler manifolds satisfying a weighted Poincar\'e inequality 
as considered in \cite{L-W4} for the Riemannian case. One then only 
requires a pointwise lower bound on the curvatures as opposed to a 
global one. Let us recall the following definition.

\proclaim{Definition E} A complete manifold  is said to have 
property ($\Cal P_\rho$) if there exists a positive function 
$\rho$ such that

$$
\int_M \rho(x)\, \phi (x) \le \int_M |\n \phi|^2(x)
$$
for all compactly supported smooth function $\phi.$  
Moreover, the conformal metric 
$$
ds_\rho^2 = \rho \, ds_M^2,
$$ 
given by multiplying the K\"ahler metric $ds_M^2$ on $M$ by $\rho$, 
is also complete.
\endproclaim

The paper is arranged as follows. In \S1, we prove Theorem A. 
In \S2, we give a systematic treatment of more general manifolds 
in the spirit of \S1. An important consequence is Theorem D, 
which deals with arbitrary irreducible locally symmetric 
spaces of noncompact type.

In \S3 and \S4, we assume that the complete K\"ahler manifold satisfies 
property ($\Cal P_\rho$).  
As pointed out in \cite{L-W1} and \cite{L-W2}, by a scaling argument, 
the pair of conditions

$$
\text{Ric}_M \ge -2(m+1) \qquad \text{and} \qquad 
\lambda_1(M) \ge \frac{m+1}2
$$
is equivalent to the pair of conditions
$$
\text{Ric}_M \ge - \frac{\lambda_1(M)}{4} \qquad \text{and} \qquad 
\lambda_1(M) >0.
$$
Written in this form, we can consider the analogue of Theorem B 
for complete K\"ahler manifolds with 
property ($\Cal P_\rho$).  In fact, Theorem 3.1 established that 
if $M$ satisfies
$$
\text{Ric}_M \ge -\frac{\rho}{4} 
$$
and if
$$
\rho(x) \to 0 \qquad \text{as} \qquad x \to \infty,
$$
then $M$ must have at most 2 nonparabolic ends.  
\footnote
{We would like to point out that recently Munteanu \cite{M}  
improved the conclusions of Theorem 3.1 and showed that $M$ 
has only one nonparabolic end.} Theorem 3.2 deals with more general
weight function $\rho$ and contains Theorem B as a special case. 
Note that when $\lambda_1(M) >0$ then an end being nonparabolic 
is equivalent to having infinite volume. So this is indeed an 
analogue of Theorem B.

In \S4, we consider the anaolgue of Theorem A for complete K\"ahler 
manifolds with property ($\Cal P_\rho$), where $\rho$ is a nonconstant
function in contrast to the situation in \S1.  
In Theorem 4.1, we show that if 
$$
BK_M \ge - \frac{\rho}{m^2} \tag 0.1
$$
and if  
$$
\rho(x) \to 0 \qquad \text{as} \qquad x \to \infty,
$$
then $M$ must have at most 2 ends providing $m \ge 3$, 
and $M$ has at most 4 ends if $m=2.$  
\footnote{Munteanu \cite{M} also improved this theorem by 
showing that $M$ has only 1 end when $m \ge 3$.} At this point, 
we would like 
to point out that the assumption on $BK_M$ can be relaxed by 
only assuming that
$$
R_{i \bar i j \bar j} \ge -1
$$
for all $i \neq j.$  In other words, we only need to assume 
a lower bound on the holomorphic bisectional curvatures 
but not the holomorphic sectional curvatures.
We also point out that we do not know 
how to deal with the case of general K\"ahler manifolds with 
property ($\Cal P_\rho$) when $\rho$ is not assumed to vanish 
at infinity. On the other hand,
we expect that as in the (real) Riemannian case \cite{L-W4}, 
it is difficult to find manifolds satisfying \thetag{0.1}
and $\rho \to \infty$ at infinity.

\heading
\S1 K\"ahler manifolds with maximum $\lambda_1$
\endheading

In this section, we concentrate on the proof of Theorem A.  
Adopting a similar notation as 
in \cite{L-W3}, we say that the K\"ahler manifold $M$ has holomorphic 
bisectional curvature bounded from below by $-C$ 
for a constant $C>0$, written as
$$
BK_M(x) \ge -C,
$$
if its curvature tensor written in any unitary frame 
$\{e_1, \dots, e_m\}$ satisfies the bound
$$
R_{i \bar i j \bar j} (x) \ge -C(1+ \delta_{ij})
$$
for all $x \in M$ and $1\le i, j\le m.$  Theorem A  can now be stated in a 
more detailed manner.

\proclaim{Theorem 1.1} Let $M^m$ be a complete K\"ahler manifold of 
complex dimension $m \ge 2$ with $\lambda_1(M) >0$.  Suppose the holomorphic 
bisectional curvature of $M$ is bounded from below by
$$
BK_M(x) \ge -\frac{\lambda_1(M)}{m^2}
$$
for all $x \in M$.  Then either
\roster
\item $M$ has only one end; or
\item $M$ is isometric to $\Bbb R \times N$ with metric
$$
ds^2_M = dt^2 + \exp(-4t)\,\omega_2^2 
+ \exp(-2t)\,\sum_{\alpha = 3}^{2m} \omega_\alpha^2,
$$
where $\{\omega_2, \dots, \omega_{2m}\}$ is an orthonormal coframe for 
a compact manifold $N.$  If $\{e_1, \dots, e_{2m}\}$ is the orthonormal frame dual to $\{dt, \omega_2, \dots, \omega_{2m}\}$ with $e_1= \frac{\p }{\p t},$ then $J e_1 = e_2.$

In the event that $M$ has bounded curvature, $N$ must be given by a compact quotient of the 
Heisenberg group.  Moreover, $M$ is isometrically covered by the complex
hyperbolic space $\Bbb CH^m.$
\endroster
\endproclaim

\demo{Proof}  By a rescaling of the metric, the assumption on the 
holomorphic bisectional curvature is equivalent to the pair of 
assumptions
$$
BK_M(x) \ge -1
$$
and
$$
\lambda_1(M) \ge m^2.
$$
For convenience sake, we will use this normalization for the purpose 
of our proof.
According to Theorem 3.2 of \cite{L-W3}, we know 
that $M$ has exactly one nonparabolic end $E_1$.  If $M$ has another 
end, $E_2$, then it must be parabolic.  Let $p \in M$ be a fixed 
point such that the compact set $B_p(R_0)$ separates the ends $E_1$ 
and $E_2,$ i.e., $E_1$ and $E_2$ are two disjoint connected 
components of $M \- B_p(R_0)$.  
Let $\gamma: [0, \infty) \rightarrow M$ be a geodesic ray 
satisfying $\gamma (0) = p$ and $\gamma (t) \to E_2(\infty),$ 
with $E_2(\infty)$ being denoted as infinity of the end $E_2.$  
We define the Busemann function $\beta$ with respect to $\gamma$ by
$$
\beta(x) = \lim_{t\to \infty} (t - r(x, \gamma(t))).
$$
The comparison theorem (Theorem 1.6) of \cite{L-W3} asserts that
$$
\d r(x, \gamma(t)) \le 2(m-1)\coth (r(x, \gamma(t))) 
+ 2 \coth (2r(x, \gamma(t))).
$$
Taking $t\to \infty,$ we conclude that
$$
\d \beta \ge -2m.
$$

For any point $x \in M$, let us consider the geodesic 
segment $\tau_t$ joining $x= \tau_t(0)$ to $\gamma(t).$  
Letting $ t \to \infty$, the sequence $\tau_t$ converges to a 
geodesic ray emanating from $x = \tau(0)$ to $ E_2(\infty).$ 
In particular, for a fixed $s >0$ and a fixed $\epsilon >0$, 
by taking sufficiently large $t$, we have
$$
r(\tau_t(s), \tau(s)) \le \epsilon.
$$
Hence, by triangle inequality,
$$
\split
\beta(\tau(s)) - \beta(\tau(0)) &=
\lim_{t\to \infty} (r(\tau(0), \gamma(t)) - r(\tau(s), \gamma(t)))\\
& = 
\lim_{t \to \infty} (r(\tau(0), \gamma(t)) - r(\tau_t(s), \gamma(t)) 
+ r(\tau_t(s), \gamma(t)) - r(\tau(s), \gamma(t)))\\
& \ge 
\lim_{t \to \infty} (r(\tau(0), \gamma(t)) - r(\tau_t(s), \gamma(t)) 
- r(\tau_t(s), \tau(s)))\\
& \ge 
\lim_{t \to \infty} (r(\tau(0), \gamma(t)) 
- r(\tau_t(s), \gamma(t))) - \epsilon\\
& = s-\epsilon.
\endsplit
$$
Since $\epsilon$ is arbitrary, we conclude that
$$
\beta(\tau(s)) - \beta(\tau(0)) \ge s.\tag1.1
$$
However, it is also clear that
$$
\split
|\beta(\tau(s)) - \beta(\tau(0))| &\le r(\tau(s), \tau(0))\\
&= s,
\endsplit
$$
hence $\beta$ is Lipschitz with Lipschitz constant 1, 
and \thetag{1.1} implies that, in fact,
$$
|\nabla \beta| = 1\tag1.2
$$
almost everywhere.  In particular, if we define the function 
$f = \exp(m \beta),$ then
$$
\split
\d f &= m \,f\, \d \beta + m^2 \, f\, |\nabla \beta|^2\\
&\ge -m^2 \, f.
\endsplit \tag 1.3
$$

Also, note that \thetag{1.1} asserts that $\beta$ when restricted 
to $\tau$ is a linear function with unit gradient.  If $x$ is in 
$M \- B_p(R_0)$ but not in $E_2$, say $x \in E_1$, then $\tau$ must 
pass through $B_p(R_0)$.  Let us denote $y$ to be the first point 
on $\tau$ that intersects $B_p(R_0)$, then \thetag{1.1} implies that
$$
\beta(y) - \beta(x) \ge r(y, x).
$$
Hence
$$
\sup_{y \in B_p(R_0)} \beta(y) - \inf_{y \in B_p(R_0)} r(y, x) 
\ge \beta(x),
$$
and combining with \thetag{1.2}, we conclude that, when restricted 
on $E_1,$ $-\beta$ is equivalent to the distance function to the 
set $B_p(R_0)$.

At this point, we would also like to point out that by tracing the 
proof of the comparison theorem in \cite{L-W3}, 
we actually proved that
$$
r_{1 \bar 1} \le \frac 12 \coth(2r)
$$
and
$$
r_{\alpha \bar \alpha} \le \frac 12 \coth (r) \qquad 
\text{for} \qquad \alpha \neq 1,
$$
where we have taken unitary frame $\{u_1, u_2, \dots, u_m\}$ with 
$$
u_1 = \frac 12\left(\nabla r - \sqrt{-1} J \nabla r\right).
$$
Hence, using the same notation as in the above discussion, 
we conclude that
$$
\beta_{i \bar i} \ge -\frac 12 \qquad 
\text{for all} \qquad 1\le i \le m,\tag1.4
$$
where $\{u_i\}$ is a unitary frame with
$$
u_1 = \frac 12 \left({\tau'(0) - \sqrt{-1} J \tau'(0)}\right).
$$
In particular, writing in terms of a real orthonormal 
frame $\{e_1, e_2, \dots, e_{2m}\}$ 
with $e_1 = \tau'(0)$, $e_2 = Je_1$, and $e_{2k}= J e_{2k-1}$ 
for $1\le k \le m,$ then, $(\beta_{ij})$, the real hessian 
of $\beta$ will satisfy
$$
\beta_{11} = 0\tag1.5
$$
since $\beta$ restricted on $\tau$ is linear with gradient 1.  
Also, \thetag{1.4} implies that
$$
\beta_{(2k-1) (2k-1)} + \beta_{(2k) (2k)} \ge -2\tag1.6
$$
for all $1 \le k \le m.$  In particular, we conclude that
$$
\beta_{22} \ge -2. \tag1.7
$$

We now claim that inequality \thetag{1.3} is indeed an equality.  
To see this, we apply the Poincar\'e inequality after multiplying 
both sides of \thetag{1.3} by $\phi^2\,f$, where $\phi$ is a 
compactly supported nonnegative cut-off function. Integrating
by parts, we conclude that
$$
\split
\lambda_1(M)\,\int_M  \phi^2\,f^2 & \le
\int_M |\n (\phi\,f)|^2\\
& = \int_M |\n \phi|^2\, f^2 + m^2 \int_M \phi^2\, f^2\\
& \qquad -\int_M \phi^2\, f (\d f+m^2\,f).
\endsplit\tag1.8
$$
We only need to justify that the first term of the right hand side 
tends to 0 for an appropriate sequence of cut-off function, then 
the hypothesis on $\lambda_1(M)$ will imply that 
$$
\int_M f\,(\d f + m^2 f) = 0.
$$ 
For $R > R_0,$ let us now choose $\phi$ to be
$$
\phi(x) = \left\{ \aligned 1\quad \quad
& \qquad \text{on} \qquad B_p(R)\\
\frac{2R-r(x)}{R}& \qquad \text{on} \qquad B_p(2R)\-B_p(R)\\
0 \quad \quad & \qquad \text{on}  \qquad M\- B_p(2R).
\endaligned \right.
$$
The first term on the right hand side of \thetag{1.8} becomes
$$
\int_M |\n \phi|^2\, f^2 = R^{-2}\int_{(B_p(2R)\-B_p(R)) \cap E_2} f^2 
+ R^{-2}\int_{(B_p(2R) \- B_p(R)) \- E_2} f^2.\tag1.9
$$
To estimate the first term, we use \thetag{1.2} for $f$ on $E_2$ 
and  Theorem 2.1 of \cite{L-W3} to obtain
$$
\split
\int_{B_p(2R)\-B_p(R) \cap E_2} f^2 
&= \sum_{i=1}^{[R]} \int_{(B_p(R+i)\-B_p(R+i-1))\cap E_2}f^2\\
&\le \sum_{i=1}^{[R]} e^{2m(R+i)}\,(V_{E_2}(R+i)\-V_{E_2}(R+i-1))\\
& \le \sum_{i=1}^{[R]} C_2\,e^{2m(R+i)}\,e^{-2m(R+i-1)}\\
&\le C_3\,R.
\endsplit\tag 1.10
$$
To estimate the second term on the right hand side of \thetag{1.9}, 
we use the fact that $-\beta$ is equivalent to the distance function 
to $B_p(R_0)$ on the other ends and Corollary 1.7 of \cite{L-W3} 
and get
$$
\split
\int_{(B_p(2R) \- B_p(R)) \- E_2} f^2 
& \le \sum_{i=1}^{[R]} \int_{(B_p(R+i) \- B_p(R+i-1))\- E_2} f^2\\
&\le \sum_{i=1}^{[R]} C_4\,e^{-2m(R+i-1)}\,e^{2m(R+i)}\\
&\le C_5\,R.
\endsplit\tag 1.11
$$
Combining \thetag{1.9}, \thetag{1.10}, and \thetag{1.11}, we 
conclude that
$$
\int_M |\n \phi|^2\,f^2 \le C_6\,R^{-1}.
$$
Letting $R \to \infty,$ we conclude our assertion and
$$
\d f =-m^2 f.
$$
In particular, all the inequalities, including \thetag{1.6} 
and \thetag{1.7}, being used to prove \thetag{1.3} are equalities 
and $f$ must be smooth.  This implies that $\beta$ is smooth 
with $|\nabla \beta| = 1$ and $\d \beta = -2m$. So $M$ is 
topologically $\Bbb R \times N$ where $N$ is diffeomorphic to the 
level set of $\beta.$  Since $M$ is assumed to have two ends, 
$N$ must be compact.  

Let us now consider the Bochner formula
$$
\d |\nabla \beta|^2 = 2 \text{Ric}_M(\nabla \beta, \nabla \beta) 
+ 2 \la \n \beta, \n \d \beta \ra + 2 \sum_{i,j}\beta_{ij}^2.
$$
Since $|\n \beta|^2 =1,$ using the assumption on the curvature 
and $\d \beta = -2m$, we conclude that
$$
0 \ge -4(m+1) + 2\sum_{i,j} \beta_{ij}^2.\tag1.12
$$
Applying \thetag{1.5}, the equality versions of \thetag{1.6} 
and \thetag{1.7}, we can estimate
$$
\split
\sum_{i,j} \beta_{ij}^2 \ge \sum_{i} \beta_{ii}^2
& \ge \beta_{22}^2 + \sum_{k+1}^{m-1} (\beta_{(2k+1)\, (2k+1)}^2 
+ \beta_{(2k+2)\, (2k+2)}^2)\\
& \ge 4 + \frac 12 \sum_{k=1}^{m-1} (\beta_{(2k +1)\, (2k +1)} 
+ \beta_{(2k+2)\, (2k+2)})^2\\
& = 4 +2(m-1)\\
& = 2(m+1).
\endsplit
$$
Hence, combining with \thetag{1.12}, we conclude that  all the 
inequalities are equalities and
the Hessian of $\beta$ is given by
$$
(\beta_{ij}) = \left( \matrix 0 & 0 & 0 & 0 & \dots & 0\\
0 & -2 & 0 & 0 & \dots & 0\\
0 & 0 & -1 & 0& \dots & 0\\
0 & 0 & 0 & -1 & \dots & 0\\
\vdots & \vdots& \vdots & \vdots &\,  &\vdots\\
0& 0& 0& 0& \dots & -1 \endmatrix \right) 
$$
Moreover, the holomorphic bisectional curvature involving 
the $u_1 = \frac 12 (e_1 - \sqrt{-1} J e_1)$ direction must be 
of the form
$$
R_{1 \bar 1 j \bar j} = -(1+ \delta_{1j}).\tag 1.13
$$

Let us now consider the level set of $\beta$ given by
$$
N_t = \{x\in M\,|\, \beta(x) = t\}.
$$
Since $|\nabla \beta| =1$, $M$ is diffeomorphic 
to $\Bbb R \times N_0$  and $e_1=\nabla \beta$ is the unit normal 
vector to $N_t$ for all $t.$  In particular, we can compute the 
second fundamental form $(h_{ij})$ of $N_t$ with respect to the 
unit normal $\n \beta$ using the hessian of $\beta$ and obtain
$$
\split
h_{ij} &= \beta_{ij} \\
&= \left\{ \aligned -2& \qquad \text{when} \qquad i = j = 2\\
-1& \qquad \text{when} \qquad i = j > 2\\
0 & \qquad \text{when} \qquad i \neq j. \endaligned \right.
\endsplit \tag 1.14
$$

Using the set of orthonormal coframe $\{\eta_1, \eta_2, \dots, \eta_{2m}\}$ 
dual to the orthonormal frame $\{e_1, e_2, \dots, e_{2m}\},$
we have the first structural equations

$$
d\eta_{i} = \eta_{ij} \wedge \eta_j,
$$
where
$\eta_{ij}$ are the connection 1-forms satisfying the condition
$$
\eta_{ij} + \eta_{ji} = 0.
$$
For $2 \le i, j \le 2m$, since 
$$
\split
\eta_{i 1}(e_j) &= \la \n_{e_j} e_i, e_1\ra\\
&= -h_{ji}
\endsplit
$$
is given by the second fundamental form of $N_t,$  we have
$$
\eta_{i 1}(e_j) = \left\{ \aligned 0
& \qquad \text{for} \qquad i \neq j\\
2& \qquad \text{for} \qquad i = j = 2\\
1& \qquad \text{for} \qquad 3 \le i= j \le 2m.
\endaligned \right.
$$
In particular, 
$$
\eta_{21} = 2 \eta_2\tag1.16
$$
and
$$
\eta_{\alpha 1} = \eta_\alpha\tag1.17
$$
for all $3\le \alpha \le 2m.$
Also note that for any vector $X$,
$$
\split
\eta_{\alpha 2}(X) &= \la \n_X e_\alpha, e_2\ra \\
& =\la \n_X e_\alpha, J e_1\ra\\
&= - \la J \n_X e_\alpha, e_1\ra\\
&= - \la \n_X Je_\alpha, e_1 \ra,
\endsplit
$$
hence using $Je_{2k+1} = e_{2k+2}$, $Je_{2k+2} = -e_{2k+1}$ 
and \thetag{1.14}, we conclude that
$$
\eta_{\alpha 2}(e_j) = \left\{ \aligned 0
& \qquad \text{for} \qquad j =1\\
-1& \qquad \text{for} \qquad \alpha = 2k+1, \, j=2k+2 
\quad \text{and for some} \quad  1\le k \le m-1\\
1& \qquad \text{for} \qquad j =2k+1,\, \alpha= 2k+2 
\quad \text{and for some} \quad 1 \le k \le m-1\\
0& \qquad \text{for} \qquad j=2. \endaligned \right.
$$
Therefore, we have
$$
\eta_{(2k+1) 2} = -\eta_{(2k+2)}\tag 1.18
$$
and
$$
\eta_{(2k+2) 2} = \eta_{(2k+1)} \tag 1.19
$$
for  $1 \le k \le m-1.$
The second structural equations assert that
$$
d\eta_{ij} - \eta_{ik} \wedge \eta_{kj} 
= \frac 12 R_{ijkl}\, \eta_l \wedge \eta_k,
$$
where $R_{ijkl}$ is the curvature tensor of $M$.  In particular, 
applying \thetag{1.16}, \thetag{1.17},\thetag{1.18}, \thetag{1.19}, 
and the first structural equation, we have
$$
\split
d\eta_{12} - \eta_{1\alpha} \wedge \eta_{\alpha 2} 
&= -2d \eta_2  
- \sum_{k=1}^{m-1} \eta_{(2k+1)} \wedge \eta_{(2k +2)} 
+ \sum_{k =1}^{m-1}\eta_{(2k+2)} \wedge \eta_{(2k+1)}\\
& = -2 \eta_{21} \wedge \eta_1 
-2\sum_{\alpha = 3}^{2m}  \eta_{2\alpha} \wedge \eta_\alpha    
+2 \sum_{k =1}^{m-1}\eta_{(2k+2)} \wedge \eta_{(2k+1)}\\
&= -4 \eta_2 \wedge \eta_1 
- 2\sum_{k=1}^{m-1} \eta_{(2k+2)} \wedge \eta_{(2k+1)}.
\endsplit
$$
This implies that
$$
R_{1212} = -4,\tag1.20
$$
$$
R_{12 (2k+1) (2k+2)} = -2\tag1.21
$$
for $1 \le k \le m-1,$ and
$$
R_{121\alpha} = 0 = R_{122 \alpha}\tag1.22
$$
for $3\le \alpha \le 2m.$  Similarly, for $1\le k \le m-1$,
$$
\split
d\eta_{1(2k+1)} &-\eta_{12} \wedge \eta_{2(2k+1)} 
- \sum_{\beta=3}^{2m}\eta_{1 \beta} \wedge \eta_{\beta (2k+1)} \\
&= - \eta_{(2k+1)} \wedge \eta_1 - \eta_{(2k+2)} \wedge \eta_2
\endsplit
$$
and
$$
\split
d\eta_{1(2k+2)} &- \eta_{12} \wedge\eta_{2(2k+2)} 
- \sum_{\beta=3}^{2m} \eta_{1 \beta} \wedge \eta_{\beta (2k+2)}\\
&= -\eta_{(2k+2)} \wedge \eta_1 + \eta_{(2k+1)} \wedge \eta_2.
\endsplit
$$

Hence,
$$
R_{1(2k+1)1(2k+1)} = -1, \tag1.23
$$
$$
R_{1(2k+1)2(2k+2)} = -1, \tag1.24
$$
and
$$
R_{1(2k+1)1i} = 0 = R_{1(2k+1)2j}
$$
for all $i \neq 2k+1$ and $j \neq 2k+2.$  We also have

$$
R_{1(2k+2)1(2k+2)} = -1, \tag1.25
$$
$$
R_{1(2k+2)2(2k+1)} = 1, \tag1.26
$$
and
$$
R_{1(2k+2)1j} = 0 = R_{1(2k+2)2i}
$$
for all $i \neq 2k+1$ and $j \neq 2k+2.$

Recall that since $|\nabla \beta|=1$ is constant along each level 
set $N_t,$ the integral curves for the vector field $\nabla \beta$ 
are all geodesics.  Moreover, the flow $\phi_t: M \rightarrow M$ 
generated by $\nabla \beta$ is a geodesic flow, and $\phi_t$ maps $N_0$
to $N_t.$  For a fixed point $p \in N_0,$ let $\tau$ be the 
geodesic given by $\tau' = \nabla \beta.$ Then the vector fields
$V_i(t)=d\phi_t(\bar{e_i})$ are Jacobi vector fields along $\tau$
for $2\le i\le 2m,$ where $\bar{e_i}$ denotes the restriction of $e_i$ 
on $N_0.$

We claim now that 
$$
V_2(t) = e^{-2t}\, e_2
$$
and for each $3 \le \alpha \le 2m,$
$$
V_\alpha(t) = e^{-t}\, e_\alpha.
$$
In particular, this implies that the metrics on the level surfaces 
$N_t$  being viewed 
as a one parameter of metrics on $N_0$ can then be written in the form

$$
ds^2_t = e^{-4t}\,\omega_2^2 
+\sum_{\alpha=3}^{2m} e^{-2t}\, \omega_\alpha^2,
$$
where $\{\omega_2, \omega_3, \omega_4, \dots, \omega_{2m}\}$ is the 
dual coframe to $\{\bar{e}_2, \bar{e}_3, \bar{e}_4, \dots, \bar{e}_{2m}\}$ 
at $N_0.$
Hence, the metric of $M$ is given by

$$
ds^2_M =dt^2+ e^{-4t}\,\omega_2^2 
+\sum_{\alpha=3}^{2m} e^{-2t}\, \omega_\alpha^2.\tag1.27
$$

Indeed, since $e_2 = Je_1 = J(\nabla \beta)$ at every point and $J$ commutes 
with the connection, $e_2$ must be a parallel vector field 
along $\tau.$  We claim that $V(t)=e^{-2t}\, e_2$
is a Jacobi field along $\tau.$  Indeed, 
$$
\n_{\tau'} (e^{-2t}\, e_2)  = -2 e^{-2t}\,e_2
$$
and
$$
\n_{\tau'} \n_{\tau'} (e^{-2t}\,e_2) = 4\,e^{-2t}\,e_2.
$$
Also, according to \thetag {1.13},
$$
\split
\la R_{\tau' e_2} \tau', e_2\ra &= R_{1212} \\
&= -4.
\endsplit
$$
Equation \thetag{1.22} also implies that
$$
\split
\la R_{\tau' e_2} \tau',e_\alpha \ra &= R_{121\alpha}\\
& = 0.
\endsplit
$$
Hence the vector field $V(t)$ satisfies the Jacobi equation
$$
\n_{\tau'} \n_{\tau'} V(t) = -R_{\tau' V} \tau'.
$$
On the other hand, 
$$
V_2(0) = d\phi_0(e_2) = e_2
$$
and
$$
\split
\n_{\tau'} V_2(0) &= \n_{V_2} e_1(p)\\
&=  \n_{e_2} e_1(p) \\
&= \sum_{j=2}^{2m}h_{2j} \,e_j \\
&= -2e_2,
\endsplit
$$
since $e_1$ and $V_2$ can be viewed as tangent vectors of a map 
from a rectangle.  Uniqueness of Jacobi field now asserts 
that $V_2(t) = V(t).$  

For each $3 \le \alpha \le 2m$, the Jacobi fields
$V_\alpha(t)$
along the geodesic $\tau$ has initial conditions
$$
V_\alpha(0) = e_\alpha\tag1.28
$$
and
$$
\split
\n_{\tau'} V_\alpha (0) &= \n_{e_\alpha} e_1(p)\\
&= \sum_{j=2}^{2m} h_{\alpha j}\, e_j \\
&= -e_\alpha.
\endsplit\tag1.29
$$
The Jacobi equation, \thetag{1.20}, and \thetag{1.22} imply that,
$$
\split
\la V_\alpha, e_2 \ra'' 
&= \la \n_{\tau'} \n_{\tau'} V_\alpha, e_2 \ra\\
&= -\la R_{e_1 V_\alpha} e_1, e_2\ra\\
&= -R_{1212}\, \la V_\alpha, e_2 \ra 
- R_{1\gamma 12} \,\la V_\alpha, e_\gamma \ra\\
&=4 \la V_\alpha, e_2\ra.
\endsplit
$$
On the other hand, using the initial conditions 
\thetag{1.28} and \thetag{1.29}, we see that
$$
\la V_\alpha, e_2 \ra(0) = 0
$$
and
$$
\split
\la V_\alpha, e_2 \ra'(0) &= -\la e_\alpha, e_2 \ra (p)\\
&= 0,
\endsplit
$$
hence we conclude that $\la V_\alpha, e_2 \ra = 0$ along $\tau.$  
In particular, $V_\alpha$ belongs to  the distribution $\Cal D$ 
spanned by the vectors $\{e_3, \dots, e_{2m}\}.$  Similarly, 
using \thetag{1.22}, \thetag{1.23}, and  \thetag{1.25} we see 
that for $3 \le \gamma \le 2m$ with $\gamma \neq \alpha$,
$$
\split
\la V_\alpha, e_\gamma\ra'' 
&= -R_{121\gamma}\,\la V_\alpha, e_2 \ra 
- R_{1\theta1\gamma} \,\la V_\alpha, e_\theta\ra\\
&=\la V_\alpha, e_\gamma \ra
\endsplit
$$
with initial conditions
$$
\la V_\alpha, e_\gamma \ra (0) = 0
$$
and
$$
\la V_\alpha, e_\gamma \ra'(0) = 0.
$$
We conclude that $\la V_\alpha, e_\gamma \ra=0$ along $\tau.$  
In particular, $V_\alpha = f(t) e_\alpha$ for some function $f(t).$  
Since the second fundamental form restricted on the subspace 
$\Cal D$ is given by the negative of the identity matrix, we 
conclude that $f(t) = e^{-t}$ and 
$$
V_\alpha(t) = e^{-t} \,e_\alpha.
$$

We will now use \thetag {1.27} to compute the curvature 
tensor of $M$ and hence $N_0.$ For the convenience sake, 
we substitute $-t$ by $t$ and rewrite \thetag {1.27} as

$$
ds^2_M =dt^2+ e^{4t}\,\omega_2^2 
+\sum_{\alpha=3}^{2m} e^{2t}\, \omega_\alpha^2.
$$

Note first the Guass curvature equation asserts that
$$
R_{ijkl} = \bar R_{ijkl} + h_{li}h_{kj} - h_{ki}h_{lj},
$$
where $\bar R_{ijkl}$ is the curvature tensor on $N_0.$
In particular, 
$$
R_{ijkl} = \left\{ \aligned  \bar R_{ijkl}
& + \delta_{li} \delta_{kj} 
- \delta_{ki} \delta_{lj} \qquad \text{if} \qquad 3\le i, j, k, l 
\le 2m\\
 \bar R_{ijkl}
& + 2 \qquad \text{if} \qquad 2=i= l \quad \text{and} 
\quad 3 \le k=j \le 2m\\
  \bar R_{ijkl}& + 2 \qquad \text{if} \qquad 2=k= j 
\quad \text{and} \quad3 \le i=l \le 2m\\
  \bar R_{ijkl}& -2 \qquad \text{if} \qquad 2=i= k
\quad \text{and} \quad 3\le j=l \le 2m\\
 \bar R_{ijkl}& -2 \qquad \text{if} \qquad 2=j= l
\quad \text{and} \quad 3\le i=k \le 2m\\
 \bar R_{ijkl}& \qquad \text{otherwise.}
\endaligned \right.\tag1.30
 $$

Since 
$$
\eta_1 = dt, 
$$
$$
\eta_2 = e^{2t}\,\omega_2, 
$$
and
$$
\eta_\alpha = e^{t}\, \omega_\alpha \qquad \text{for} \qquad \alpha 
= 3, \dots 2m,
$$
we obtain
$$
d\eta_1 = 0,\tag 1.31
$$
$$
\split
d \eta_2 &= 2e^{2t}\,\eta_1\wedge \omega_2 
+ e^{2t}\,\sum_{\alpha=3}^{2m}\omega_{2\alpha} \wedge \omega_\alpha\\
&= -2 \eta_2 \wedge \eta_1 
+ e^{t}\,\sum_{\alpha=3}^{2m} \omega_{2\alpha} \wedge \eta_\alpha,
\endsplit\tag 1.32
$$
and
$$
\split
d \eta_\alpha &= e^{t}\, \eta_1 \wedge \omega_\alpha 
+ e^{t} \, \omega_{\alpha 2}\wedge \omega_2 
+ e^{t}\, \omega_{\alpha \beta} \wedge \omega_\beta\\
&= -\eta_\alpha \wedge \eta_1 
+ e^{-t}\, \omega_{\alpha 2} \wedge \eta_2 
+ \omega_{\alpha \beta} \wedge \eta_\beta,
\endsplit\tag 1.33
$$
where $\omega_{ij}$ are the connection forms of $N_0.$  In the 
above and all subsequent computations, we will adopt the convention 
that $3 \le \alpha, \beta \le 2m,$ $2 \le i, j \le 2m,$ 
$ 2\le s, t \le m,$ and $1\le p, q \le 2m.$  

Note that by \thetag {1.18} and \thetag {1.19},
$$
\omega_{2(2s-1)} = -\omega_{(2s)}\tag 1.34
$$
and
$$
\omega_{2(2s)} = \omega_{(2s-1)}. \tag 1.35
$$

Equation \thetag{1.32} and \thetag{1.33} imply that the connection 
forms are given by
$$
\split
\eta_{12} &=-\eta_{21}\\
&= 2 \eta_2,
\endsplit \tag 1.36
$$
$$
\split
\eta_{1\alpha} &=-\eta_{\alpha 1}\\
&= \eta_\alpha,
\endsplit \tag 1.37
$$
$$
\split
\eta_{(2s-1) \beta} &= -\eta_{\beta (2s-1)}\\
&= \left\{ \aligned \omega_{(2s-1) \beta} +(1-&e^{-2t})\,
\eta_2 \qquad \text{if} \qquad \beta= 2s\\
\omega_{(2s-1) \beta}& \qquad \text{if} \qquad \beta \neq 2s, 
\endaligned \right.
\endsplit \tag 1.38
$$
$$
\split
\eta_{(2s) \beta} &= -\eta_{\beta (2s)}\\
&= \left\{ \aligned \omega_{(2s) \beta} 
-(1-&e^{-2t})\,\eta_2 \qquad \text{if} \qquad \beta= 2s-1\\
\omega_{(2s) \beta}& \qquad \text{if} 
\qquad \beta \neq 2s-1, \endaligned \right.
\endsplit \tag 1.39
$$
and
$$
\split
\eta_{2 \alpha} &=-\eta_{\alpha 2}\\
&=e^t\, \omega_{2\alpha}.
\endsplit \tag 1.40
$$
Indeed, if we substitute \thetag{1.36} and \thetag{1.40} into 
the first structural equation
$$
d \eta_2 =\eta_{21}\wedge \eta_1 +  \eta_{2\alpha} \wedge \eta_{\alpha}
$$
we obtain \thetag{1.32}.
Also, using \thetag{1.34}, \thetag{1.35}, \thetag{1.37}, 
\thetag{1.38},and \thetag{1.40}, we have
$$
\split
d\eta_{(2s-1)} &= \eta_{(2s-1) 1} \wedge \eta_1 
+\eta_{(2s-1) 2} \wedge \eta_2 
+ \eta_{(2s-1) \beta} \wedge \eta_\beta\\
&=- \eta_{(2s-1)} \wedge \eta_1 - e^t\,\omega_{2(2s-1)}\wedge \eta_2 \\
& \qquad + \omega_{(2s-1) \beta} \wedge \eta_\beta 
+ (1-e^{-2t}) \,\eta_2 \wedge \eta_{(2s)} \\
&= - \eta_{(2s-1)} \wedge \eta_1 +e^{-2t}\, \eta_{(2s)} \wedge \eta_2 
+ \omega_{(2s-1) \beta} \wedge \eta_\beta,
\endsplit
$$
validating \thetag{1.33}. A similar computation also 
validates \thetag{1.39}

To compute the curvature, we consider the second structural 
equations.  In particular,
$$
\split
d \eta_{12} - \eta_{1 \alpha} \wedge \eta_{\alpha 2} 
&= 2d\eta_2 - \sum_{s=2}^m\eta_{(2s-1)} \wedge   \eta_{(2s)}  
+\sum_{s=2}^m \eta_{(2s)} \wedge \eta_{(2s-1)}\\
& = -4 \eta_2 \wedge \eta_1  
+2\sum_{s=2}^m \eta_{(2s-1)} \wedge \eta_{(2s)},
\endsplit
$$
hence
$$
R_{12ij}=\left\{\aligned  -4& \qquad \text{if} \qquad i=1, j=2\\
-2& \qquad \text{if} \qquad i=2s-1, j=2s\\
0& \qquad \text{otherwise.}
\endaligned \right.\tag{1.41}
$$
Also,
$$
\split
d \eta_{1\alpha} - \eta_{12} \wedge \eta_{2 \alpha} 
- \eta_{1 \beta} \wedge \eta_{\beta \alpha} &=
\eta_{\alpha 1} \wedge \eta_1 + \eta_{\alpha 2} \wedge \eta_2 
+ \eta_{\alpha \beta} \wedge \eta_\beta 
- 2\eta_{2} \wedge \eta_{2 \alpha} 
- \eta_{\beta} \wedge \eta_{\beta \alpha}\\
&= -\eta_{\alpha} \wedge \eta_1  
-e^t\, \eta_2 \wedge  \omega_{2\alpha}\\
&= \left\{ \aligned -\eta_{(2s-1)} \wedge \eta_1 
&+ \eta_2 \wedge \eta_{(2s)}\qquad \text{if} \qquad \alpha = 2s-1\\
-\eta_{(2s)} \wedge \eta_1 & - \eta_2 \wedge \eta_{(2s-1)} 
\qquad \text{if} \qquad \alpha = 2s
\endaligned \right.
\endsplit
$$
hence
$$
R_{1 \alpha ij}= \left\{ \aligned -1& \qquad \text{if} \qquad i=1, j=\alpha\\
1& \qquad \text{if} \qquad \alpha=2s-1, i=2s, j=2\\
-1&\qquad \text{if} \qquad \alpha=2s, i=2s-1, j=2\\
0& \qquad \text{otherwise.} \endaligned \right.\tag1.42
$$

Moreover,
$$
\split
d\eta_{2\alpha} &- \eta_{21} \wedge \eta_{1 \alpha} 
- \eta_{2\beta} \wedge \eta_{\beta \alpha}\\
&= \left\{ \aligned d({e^t}\, \omega_{2\alpha}) 
&+2\eta_2 \wedge \eta_\alpha 
- e^t\, \omega_{2\beta}\wedge \omega_{\beta \alpha} 
+e^t (1-e^{-2t})\, \omega_{2(2s)}\wedge  \eta_2, \qquad \alpha= 2s-1\\
d({e^t}\, \omega_{2\alpha}) &+2\eta_2 \wedge \eta_\alpha 
- e^t\, \omega_{2\beta}\wedge \omega_{\beta \alpha} 
- e^t(1-e^{-2t})\, \omega_{2(2s-1)}\wedge  \eta_2, \qquad \alpha= 2s
\endaligned \right.\\
&=  e^t \,\omega_1\wedge \omega_{2\alpha} 
+ \frac 12 e^t \bar R_{2\alpha ij}\,\omega_j\wedge \omega_i 
-2 \eta_\alpha \wedge \eta_2 
+ e^t(1-e^{-2t})\,\omega_{\alpha}\wedge \eta_2\\
&= \left\{ \aligned -\eta_1 \wedge \eta_{(2s)} 
&+ \frac 12 e^t \bar R_{2\alpha ij}\,\omega_j\wedge \omega_i 
-2 \eta_\alpha \wedge \eta_2 
+ (1-e^{-2t})\,\eta_{\alpha}\wedge \eta_2, \qquad \alpha = 2s-1\\
 \eta_1 \wedge \eta_{(2s-1)}
& +\frac 12 e^t \bar R_{2\alpha ij}\,\omega_j\wedge \omega_i 
-2 \eta_\alpha \wedge \eta_2 
+ (1-e^{-2t})\,\eta_{\alpha}\wedge \eta_2, \qquad \alpha = 2s, 
\endaligned \right.
\endsplit\tag 1.43
$$
where $\bar R_{23ij}$ is the curvature tensor of $N_0.$
In particular,
$$
K_M(e_2, e_\alpha) = e^{-2t}\,K_0(e_2, e_\alpha) -1 - e^{-2t},\tag 1.44
$$
where $K_0$ is the sectional curvature of $N_0$, On the other 
hand, \thetag{1.41} and \thetag{1.42} together with the K\"ahler condition imply that
the curvature tensor involving the $e_2$ direction is completely determined.
$$
\split
K_M(e_2, e_\alpha) &= K_M(e_1, Je_{\alpha})\\
=-1.
\endsplit 
$$
Combining with \thetag{1.44}, we conclude that
$$
K_0(e_2, e_\alpha) = 1
$$
for $\alpha=3, \dots, 2m.$  Equation \thetag{1.43} also implies that
$$
R_{2\alpha 2\beta} = e^{-2t}\, \bar R_{2\alpha 2\beta} 
\qquad \text{for} \qquad \beta \neq \alpha
$$
and
$$
R_{2\alpha \beta \gamma} 
= e^{-t}\,\bar R_{2\alpha \beta \gamma}.
$$
Again, using
$$
R_{2 (2s) ij} = R_{1(2s-1)ij}
$$
$$
R_{2(2s-1)ij}=-R_{1(2s)ij}
$$
and \thetag{1.42}, we conclude that
$$
R_{2\alpha 2 \beta}=0 \qquad \text{for} \qquad \beta \neq \alpha
$$
and
$$
R_{2\alpha \beta \gamma} = 0,
$$
hence
$$
\bar R_{2\alpha 2\beta}= 0 \qquad \text{for} \qquad \beta \neq \alpha
$$
and
$$
\bar R_{2\alpha \beta \gamma}=0.
$$

It remains for us to compute the curvature tensor in the directions involving only $e_\alpha$ for $3 \le \alpha \le 2m.$
Following the computation of the second structural equations, 
using \thetag {1.34}, \thetag {1.36}, \thetag {1.37}, 
and \thetag {1.40}, we have
$$
\split
d\eta_{(2s-1) (2s)} &- \eta_{(2s-1) 1}\wedge \eta_{1(2s)} 
- \eta_{(2s-1) 2} \wedge \eta_{2 (2s)} 
- \eta_{(2s-1) \gamma} \wedge \eta_{\gamma (2s)}\\
&=d\omega_{(2s-1)(2s)} + 2e^{-2t}\, \eta_1 \wedge \eta_2  
+ (1-e^{-2t}) \,  d\eta_2 + \eta_{(2s-1)} \wedge \eta_{(2s)} \\
&\qquad + e^{2t}\, \omega_{2(2s-1)}\wedge \omega_{2(2s)} 
-\omega_{(2s-1)\gamma} \wedge \omega_{\gamma (2s)}\\
&= d\omega_{(2s-1)(2s)} -\omega_{(2s-1)2} \wedge \omega_{2(2s)}
-\omega_{(2s-1)\gamma} \wedge \omega_{\gamma (2s)} 
+ 2e^{-2t}\, \eta_1 \wedge \eta_2  \\
&\qquad -2 (1-e^{-2t}) \,  \eta_2\wedge \eta_1
+(1-e^{-2t}) e^t \,\omega_{2\gamma} \wedge \eta_\gamma\\ 
&\qquad + \eta_{(2s-1)} \wedge \eta_{(2s)}
+(1 - e^{2t})\, \omega_{(2s)}\wedge \omega_{(2s-1)}\\
&= \frac 12 \bar R_{(2s-1)(2s)ij}\, \omega_j \wedge \omega_i 
+ e^{-2t}\, \eta_{(2s)} \wedge \eta_{(2s-1)} + 2\eta_1 \wedge \eta_2\\
&\qquad - 2(1-e^{-2t})\, \sum_{r=2}^m\eta_{(2r)} \wedge \eta_{(2r-1)} 
- 2\eta_{(2s)} \wedge \eta_{(2s-1)}.
\endsplit\tag1.45
$$
This implies that
$$
K_M(e_{2s-1}, e_{2s}) = e^{-2t}\,(K_0(e_{2s-1}, e_{2s}) +3) -4.\tag 1.46
$$
Equation \thetag{1.45} also implies that
$$
R_{(2s-1)(2s)(2r-1)(2r)} 
= e^{-2t}\,(\bar R_{(2s-1)(2s)(2r-1)(2r)} +2) -2\tag 1.47
$$
for $r \neq s,$ and
$$
R_{(2s-1)(2s)\alpha \beta} 
= e^{-2t} \, \bar R_{(2s-1)(2s)\alpha \beta}\tag1.48
$$
for $\alpha \neq 2s-1, 2s$, $\beta \neq 2s-1, 2s$, 
and $\alpha \neq 2r-1$ when $\beta =2r.$ 
For $r \neq s$, we compute
$$
\split
d\eta_{(2s-1)(2r-1)}& -\eta_{(2s-1)1} \wedge \eta_{1(2r-1)} 
- \eta_{(2s-1)2} \wedge \eta_{2(2r-1)} 
- \eta_{(2s-1) \gamma} \wedge \eta_{\gamma (2r-1)}\\
&= d\omega_{(2s-1)(2r-1)}  +\eta_{(2s-1)} \wedge \eta_{(2r-1)} 
+ \eta_{(2s)} \wedge \eta_{(2r)} 
- \omega_{(2s-1)\gamma} \wedge \omega_{\gamma (2r-1)} \\
&\qquad - (1-e^{-2t})\, \eta_2 \wedge \omega_{(2s)(2r-1)} 
+ \omega_{(2s-1) (2r)} \wedge (1-e^{-2t})\, \eta_2\\
&= \frac 12 \bar R_{(2s-1)(2r-1)ij}\, \omega_j \wedge\omega_i 
- \omega_{(2s)}\wedge \omega_{(2r)} + \eta_{(2s-1)} \wedge \eta_{(2r-1)} 
+ \eta_{(2s)} \wedge \eta_{(2r)}\\
&= \frac 12 \bar R_{(2s-1)(2r-1)ij}\, \omega_j \wedge\omega_i 
- e^{-2t}\,\eta_{(2s)}\wedge \eta_{(2r)} + \eta_{(2s-1)} \wedge \eta_{(2r-1)} 
+ \eta_{(2s)} \wedge \eta_{(2r)},
\endsplit 
$$
where we have used the fact that the K\"ahler condition implies that
$$
\omega_{(2s)(2r-1)} = -\omega_{(2s-1)(2r)}.
$$
This implies that
$$
K_M(e_{2s-1}, e_{2r-1})= e^{-2t}\,K_0(e_{2s-1}, e_{2r-1}) -1,\tag4.49
$$
$$
R_{(2s-1)(2r-1)(2s)(2r)} 
= e^{-2t}\,(\bar R_{(2s-1)(2r-1)(2s)(2r)} +1) -1. \tag4.50
$$
A similar computation also yields
$$
\split
d\eta_{(2s)(2r)} &-\eta_{(2s)1} \wedge \eta_{1(2r)} 
- \eta_{(2s)2}\wedge \eta_{2(2r)} 
-\eta_{(2s)\gamma} \wedge \eta_{\gamma (2r)}\\
&= \frac12 \bar R_{(2s)(2r)ij}\, \omega_j \wedge \omega_i 
+ (e^{-2t}-1)\, \eta_{(2s-1)}\wedge \eta_{(2r-1)} - \eta_{(2r)}\wedge \eta_{(2s)},
\endsplit
$$
implying that
$$
K_M(e_{2s}, e_{2r}) = e^{-2t}\,K_0(e_{2s}, e_{2r}) -1.\tag4.51
$$
Finally, we compute for $s\neq r,$ 
$$
\split
d\eta_{(2s-1)(2r)} &-\eta_{(2s-1)1} \wedge \eta_{1(2r)} 
- \eta_{(2s-1)2}\wedge \eta_{2(2r)} 
-\eta_{(2s-1)\gamma} \wedge \eta_{\gamma (2r)}\\
&= \frac12 \bar R_{(2s-1)(2r)ij}\, \omega_j \wedge \omega_i 
+ \eta_{(2s-1)}\wedge \eta_{(2r)} - \eta_{(2s)}\wedge \eta_{(2r-1)},
\endsplit
$$
implying that
$$
K_M(e_{2s-1}, e_{2r}) = e^{-2t}\,K_0(e_{2s-1}, e_{2r}) -1.\tag4.52
$$

If $M$ has bounded curvature, then equations \thetag{4.16}, \thetag{4.17}, \thetag{4.18}, \thetag{4.19}, \thetag{4.50}, \thetag{4.51}, 
and \thetag{4.52}  assert that all the terms involving the $e^{-2t}$ factor must be zero.  Hence, we conclude that
$$
K_M(e_{2s-1}, e_{2s})= -4,
$$
$$
R_{(2s-1)(2s)(2r-1)(2r)} =  -2 \qquad \text{for} \qquad r \neq s,
$$
$$
R_{(2s-1)(2s)\alpha \beta} 
=0
$$
for $\alpha \neq 2s-1, 2s$, $\beta \neq 2s-1, 2s$, 
and $\alpha \neq 2r-1$ when $\beta =2r.$
Also,
$$
K_M(e_{2s-1}, e_{2r-1})=  -1,
$$
$$
R_{(2s-1)(2r-1)(2s)(2r)} 
=  -1,
$$
$$
K_M(e_{2s}, e_{2r}) = -1,
$$
and
$$
K_M(e_{2s-1}, e_{2r}) =-1.
$$
In particular, these determined the whole curvature tensor for 
$M$ and $N_0,$  and $M$ must have constant holomorphic 
bisectional curvature, hence must be covered by $\Bbb {CH}^m.$
In this case, to see that $N$ is a compact quotient of the Heisenberg group, 
one first observes that since $\beta$ has no critical point $N$ 
must be a compact quotient of a horosphere of $\Bbb {CH}^m.$  
It is then not difficult to see (see \cite{B-DR}) that a horosphere 
is given by the Heisenberg group. 
\qed
\enddemo
 
We should take this opportunity to point out that since the lattice 
consisting of even integers in $\Bbb R^{2m-1}$ is a discrete 
subgroup of the Heisenberg group and their quotient is obviously 
compact, this gives an example of the existence of case (2) in the 
conclusion of Theorem 1.1 when $M$ has bounded curvature.  

One can also construct an example of case (2) when $M$ has unbounded curvature.  
Indeed, let us consider $N=S^{2m-1}$ the unit sphere in $\Bbb C^m$ with the
induced contact $1$-form $\omega_2,$ that is, $\omega_2=J dr,$ where $J$ is the
standard complex structure on $\Bbb C^m.$ Let $\{\omega_2, \dots \omega_{2m}\}$ be 
an orthonormal coframe of $N$ such that $J \omega_{2s-1} = \omega_{2s}$ for $2\le s\le m.$
Since the K\"ahler form on $\Bbb C^m$ is given by

$$
\omega=r\,dr \wedge \omega_2+r^2 (\omega_3 \wedge \omega_4+\dots+\omega_{2m-1} \wedge \omega_{2m})
$$
and $d \omega=0,$ one concludes that

$$
d\omega_2=2(\omega_3 \wedge \omega_4+\dots+\omega_{2m-1} \wedge \omega_{2m}) \tag 4.53
$$
on $N=S^{2m-1}.$
Now, we consider a metric on $M = \Bbb R \times N$ given by 
$$
ds_M^2 = dt^2 + e^{4t}\,\omega_2^2 + e^{2t}\,\sum_{\alpha=3}^{2m} \omega_\alpha^2 \tag 4.54
$$
with the almost complex structure defined by $J dt=\omega_2$ and $J \omega_{2s-1} = \omega_{2s}$ for $2\le s\le m.$  
One checks readily that this almost complex structure is integrable. Also, using \thetag{4.53},
one concludes by direct computation that the K\"ahler form associated to the metric \thetag{4.54}, which is given by

$$
\omega_M=e^{2t}dt\wedge \omega_2+e^{2t}(\omega_3 \wedge \omega_4+\dots+\omega_{2m-1} \wedge \omega_{2m}),
$$ 
must be closed.
Hence, $M$ is a K\"ahler manifold. Finally, from the curvature computations carried out above, one sees that 
it satisfies all the conditions of Theorem 1.1. 

We remark that this construction works for any compact hypersurface $N$ of a K\"ahler manifold so long as
the induced contact structure on $N$ satisfies \thetag{4.53} and $N$ also satisfies a suitable curvature lower bound.
 
\heading
\S2 Locally Symmetric Spaces
\endheading

The argument in \S1 can be generalized to the following situation.

\proclaim{Theorem 2.1} Let $M^n$ be a complete Riemannian manifold 
of dimension $n$.  Suppose  $f:(0,\infty) \rightarrow \Bbb R$ is a 
function with the property that 
$$
\lim_{r\to \infty} f(r) = 2a>0
$$
and
$$
\int_{R_0}^\infty (f(r) - 2a)\, dr < \infty
$$
for some $R_0< \infty$.  Assume that for any point $p\in M$, 
and if $r(x)$ is the distance function to the point $p$, we have
$$
\d r(x) \le f(r(x))
$$
in the weak sense.  If $M$ has at least one parabolic end, then
$$
\lambda_1(M) \le a^2.
$$
Moreover, if $\lambda_1(M) = a^2,$ then let $\gamma(t)$ be a 
geodesic ray issuing from a fixed point $p$ to infinity of the 
parabolic end, and the Buseman function 
$$
\beta(x) = \lim_{t \to \infty} (t - r(\gamma(t), x))
$$
with respect to $\gamma$ must satisfy
$$
\d \beta = a^2,
$$
$$
|\nabla \beta|=1.
$$
Hence, $M$ must be homeomorphic to $\Bbb R \times N$ for some 
compact manifold $N$ given by the level set of $\beta.$
\endproclaim

\demo{Proof} We may assume that $\lambda_1(M) >0$ as otherwise the 
theorem is trivial.  In this case,  $M$ is nonparabolic.  
By assumption that $M$ has at least one parabolic end, let us 
denote $E_2$ to be a parabolic end. Then $E_1=M\-E_2$ is a 
nonparabolic end.  Following the argument of Theorem 1.1, 
let $\gamma:[0, \infty) \rightarrow M$ be a geodesic ray 
with $\gamma(0)= p$ and $\gamma(t) \to E_2(\infty).$  Also, 
using the inequality on $\d r(x),$ we conclude that
$$
\d \beta \ge -2a\tag{2.1}
$$
and
$\beta$ is Lipschitz with Lipschitz constant 1.  
Setting $f = \exp (a\beta)$ and using \thetag{2.1}, a direct 
computation yields
$$
\d f \ge - a^2 \,f.
$$
Following the same argument as in the proof of Theorem 1.1, we obtain
$$
(\lambda_1(M) - a^2)\int_M \phi^2\,f^2 
\le \int_M |\n \phi|^2\,f^2.\tag2.2
$$
Assuming the contrary that $\lambda_1 > a^2$, we obtain a 
contradiction if we can justify the right hand side tends to 0 
for a sequence of cut-off functions $\phi$ unless 
$\lambda_1(M) = a^2$ and all the above inequalities are equalities.  

As in the proof of Theorem 1.1, to estimate the right hand side 
of \thetag{2.2} on the parabolic end $E_2$, it suffices to show that
$$
V_{E_2}(R) \- V_{E_2}(R-1) \le \exp(-2a (R-1)).
$$
This follows from the volume estimate (Theorem 2.1 of \cite{L-W3})
$$
\split
V_{E_2}(R) \- V_{E_2}(R-1) &\le \exp (-2\sqrt{\lambda_1(M)}(R-1))\\
&\le \exp(-2a(R-1)).
\endsplit
$$
For the non-parabolic end $E_1$, we need the fact that
$$
B_p(R)\-B_p(R-1) \le \exp(2aR).
$$
To see this, using  the fact that $\d r$ is the mean curvature 
$H(\theta, r)$ of $\p B_p(r)$ in terms of polar coordinates 
centered at $p,$ we have
$$
J(\theta, r) 
\le J(\theta, r_0)\, \exp\left(\int_{r_0}^r f(t)\,dt\right), 
$$
where $J$ is the area element in polar coordinates. 
Integrating over the variable $\theta$, we have

$$
A_p(r) \le A_p(R_0) \,\exp\left(\int_{R_0}^r f(t)\, dt\right),
$$
where $A_p(r)$ denotes the area of the set $\p B_p(R).$
Hence using the assumption on $f$, we conclude that
$$
\split
(V_p(R)-V_p(R-1))&\le 
A_p(R_0)\,\int_{R-1}^R \exp(\int_{R_0}^r f(t)\,dt)\,dr \\
&\le A_p(R_0) \,\int_{R-1}^R \exp(2a(r-R_0) + C)\,dr\\
&\le C_1\,A_p(R_0) \,\exp(2aR)
\endsplit
$$
as needed.  In conclusion, $\lambda_1(M) \le a^2.$

Following the proof of Theorem 1.1, if 
$$
\lambda_1(M) = a^2,
$$
then we conclude that
$$
\d \beta = -2a,\tag 2.3
$$
$$
|\n \beta|=1,\tag 2.4
$$
and
$\beta$ has no critical points.  In particular, $M$ must be 
homeomorphic to $\Bbb R \times N$ for some compact manifold $N.$

The Bochner formula together with \thetag{2.3} and \thetag{2.4} 
implies that
$$
\split
0 &= \d |\n \beta|^2\\
&= 2 \beta_{ij}^2 + 2\la \n \beta, \n \d \beta\ra 
+ 2 \text{Ric}_{11} \\
&= 2 \beta_{ij}^2 + 2 \text{Ric}_{11}
\endsplit
$$
for the unit vector $e_1 = \n \beta$.  Using \thetag{2.4} again, 
this implies that $\beta_{1i}=0$ for all $i$ and the second 
fundamental form $II$ of a level set of $\beta$ satisfies
$$
|II|^2= -\text{Ric}_{11}.\tag2.5
$$
On the other hand, the Gauss curvature equation asserts that
for $\alpha \neq \tau,$ 
$$
K_N(e_\alpha, e_\tau) = K_M(e_\alpha, e_\tau)
+\lambda_{\alpha}\lambda_{\tau}
$$
for an orthonormal frame $\{e_\alpha\}_{\alpha=2}^n$ on the level 
set of $\beta$ that diagonalizes $II$ with corresponding eigenvalues
$\{\lambda_\alpha\}_{\alpha=2}^n.$
Since the scalar curvature of $M$ is given by
$$
\split
S_M
&=\sum_{i=1}^n \text{Ric}_{ii}\\
&= \text{Ric}_{11} + \sum_{\alpha=2}^n \text{Ric}_{\alpha \alpha}\\
&= 2\text{Ric}_{11} + \sum_{\alpha,\,
\tau \neq 1} K_M(e_\alpha,  e_\tau)\\
&=2\text{Ric}_{11}  +\sum_{\tau \neq \alpha}
K_N(e_\alpha,  e_\tau) 
- \sum_{\tau \neq \alpha}\lambda_\alpha \lambda_\tau.\\
&= 2\text{Ric}_{11}  + S_N 
- \sum_{\tau \neq \alpha}\lambda_\alpha \lambda_\tau,
\endsplit
$$
this implies that
$$
S_N -S_M +2\text{Ric}_{11} 
=  \sum_{\tau \neq \alpha}\lambda_\alpha \lambda_\tau.\tag2.6
$$
On the other hand, \thetag{2.3} and \thetag{2.4} 
assert that
$$
H= -2a
$$
where $H$ is the mean curvature of the level set of $\beta.$  
Combining with \thetag{2.5} and \thetag{2.6}, we conclude that
$$
\split
4a^2 &= H^2\\
&= |II|^2 + \sum_{\tau \neq \alpha} \lambda_{\alpha}\lambda_{\tau}\\
&= S_N - S_M +\text{Ric}_{11}.
\endsplit
$$
Hence
$$
S_N = 4a^2 +S_M -\text{Ric}_{11}\tag2.7
$$
and
$$
\sum_{\tau \neq \alpha} \lambda_{\alpha}\lambda_{\tau}
= 4a^2 + \text{Ric}_{11}.
$$
Also note that the inequality 
$$
|II|^2 \ge \frac{H^2}{n-1}
$$
implies that
$$
-(n-1)\text{Ric}_{11} \ge 4a^2\tag2.8
$$
with equality if and only if $\lambda_{\alpha} = \lambda_{\tau}$ 
for all $\alpha$ and $\tau.$  \qed
\enddemo

We first observe that the above theorem allows us to recover a 
theorem proved in \cite{L-W2}.
\proclaim{Corollary 2.2}  Let $M^n$ be a complete manifold of 
dimension $n \ge2.$  Assume that
$$
\text{Ric}_M \ge -(n-1)
$$
and
$$
\lambda_1(M) \ge \frac{(n-1)^2}{4}.
$$
Then
$M$ must either have no finite volume end or it must be a warped 
product $M= \Bbb R \times N$ with metric given by
$$
ds_M^2 = dt^2 + exp(2t)\,ds_N^2,
$$
where $N$ is a compact manifold whose metric $ds_N^2$ has 
nonnegative Ricci curvature.
\endproclaim

\demo{Proof} We first observe that the assumption on the Ricci 
curvature and Laplacian comparison theorem asserts that
$$
\d r \le (n-1)\coth r,
$$
hence one checks readily that the function $f(r)= (n-1)\coth r$ 
satisfies the hypothesis of Theorem 2.1 with $a=\frac{n-1}2.$  
Therefore we conclude that if $M$ has a parabolic end it must be 
homeomorphic to $\Bbb R \times N$ for some compact manifold.  
Moreover, since $\lambda_1(M) = \frac{(n-1)^2}4 >0$, an end 
being parabolic is equivalent to having finite volume.  
Also,  \thetag{2.8} takes the form
$$
\split
-(n-1) &\ge \text{Ric}_{11}\\
&\ge -(n-1).
\endsplit
$$
In particular, \thetag{2.8} becomes an equality and we conclude that
$$
II = -\left(\delta_{\alpha \tau}\right)
$$
is a diagonal matrix.  In this case, the metric on $M$ must be of 
the form
$$
ds_M^2= dt^2 + \exp(-2t)\,ds_N^2.
$$
A direct computation shows that the sectional curvatures for the 
sections span by $\{e_1= \frac{\p }{\p t}, e_{\alpha}\}$ is given by
$$
K_M(e_1, e_\alpha) = -1.
$$
The Gauss curvature equation implies that
$$
K_M(e_\alpha , e_\tau)=K_N(e_\alpha, e_\tau) - 1,
$$
and hence
$$
\text{Ric}_{\alpha \alpha} 
=\overline{\text{Ric}}_{\alpha \alpha} -(n-1),
$$
where $\overline{\text{Ric}}_{\alpha \alpha}$ is the Ricci 
curvature of $N$.  This implies that
$$
\overline{\text{Ric}}_{\alpha \alpha} \ge 0.
$$
The theorem follows by setting $t$ to be $-t.$ \qed
\enddemo

Let us now focus our discussion on a special case of Theorem 2.1 
when $M$ is an Einstein manifold with Einstein constant $-\gamma<0.$  
In particular,  \thetag{2.7} and \thetag{2.8} become
$$
S_N= 4a^2 - (n-1)\gamma
$$
and
$$
(n-1)\gamma \ge 4a^2.
$$
This implies that $S_N \le 0$ with $S_N=0$ if and only if 
$$
II = -\frac{2a}{n-1} \delta_{\alpha \tau}.
$$
In the case of equality, using the same argument as in the above corollary, 
we conclude that 
$$
ds^2_M = dt^2 + \exp\left( -\frac{4a}{n-1}t\right) \,ds_N^2,
$$
hence
$$
K(e_1, e_\alpha)=-\frac{4a^2}{(n-1)^2}.
$$
and
$$
\split
-\frac{4a^2}{n-1} &=-\gamma \\
&= \text{Ric}_{\alpha \alpha}\\
&= \overline{\text{Ric}}_{\alpha \alpha} - \frac{4a^2}{n-1}.
\endsplit
$$
We conclude that $N$ must be Ricci flat.

The next theorem gives a more detailed description on the conclusion of Theorem D.

\proclaim{Theorem 2.3} Let $M$ be an irreducible locally symmetric 
space covered by an irreducible symmetric space of noncompact 
type $\tilde M = G/K.$  Suppose $\lambda_1(M) = \lambda_1(\tilde M).$  
Then either 
\roster
\item $M$ has only one end; or
\item  $M$ is isometric to $\Bbb R \times N$ with metric 
$$
ds^2_M = dt^2 + \sum_{\alpha=2}^m \exp\left(-2b_\alpha t\right)\,
\omega_\alpha^2,
$$
where $\{\omega_2, \dots \omega_n\}$ is an orthonormal basis for a 
compact manifold $N$ given by a compact quotient of the horosphere 
of $\tilde M$ and $b_\alpha$ are the nonnegative constants
such that $\{-b_\alpha^2\}$ are the eigenvalues of the symmetric tensor
$$
A_{\alpha \gamma} = R_{1\alpha 1\gamma}
$$
for a fixed direction $e_1.$
\endroster
\endproclaim

\demo{Proof} If $\tilde M$ is of rank one, then it must be either 
the real hyperbolic space, the complex hyperbolic space, the 
quaternionic hyperbolic space, or the Cayley plane.  For the case 
when $\tilde M$ is the real hyperbolic space, the theorem follows 
from the previous work of the authors \cite{L-W1} and \cite{L-W2} 
as indicated by Corollary 2.2.  In this case, the cross section is 
a flat manifold since $M$ is assumed to have constant $-1$ curvature 
and the horosphere of $\Bbb H^n$ are simply Euclidean space 
$\Bbb R^{n-1}.$  In the case when $\tilde M$ is the complex 
hyperbolic space, this was covered by Theorem 1.1.  The remaining 
two rank one cases given by the quaternionic hyperbolic space and 
the Cayley plane are separately studied in \cite{K-L-Z} and \cite{Lm}.

We may now concentrate our attention on the cases when $\tilde M$ 
is an irreducible symmetric space with rank at least 2.  In this case, 
we first observe that using a formula of Matsushima \cite{Ma} (generalized by 
Jost and Yau \cite{J-Y} to harmonic maps) we can rule out the 
existence of a second nonparabolic end for $M.$  Indeed, if $M$ 
has two nonparabolic ends, then the Tam-Li theory asserts the 
existence of a nonconstant bounded harmonic function, $f$, with 
finite Dirichlet integral.   A special case of Lemma 1 and Theorem 2 
in \cite{J-Y} now implies that if $M$ is a compact irreducible 
locally symmetric space of noncompact type with rank at least 2, 
then the hessian of $f$ must be identically 0.  To apply their 
argument to our situation, we only needs to justify integration 
by parts in certain steps of their proof.  However, as in previous 
argument in this paper, it is sufficient to check that the norms of 
the gradient $|\nabla f|$ and the hessian $|f_{ij}|$ 
are both in $L^2.$  The first follows from the fact that $f$ has 
finite Dirichlet integral and the second follows from the standard 
cut-off argument applied to the Bochner formula
$$
\d |\nabla f|^2 \ge 2 \text{Ric}_{ij}\, f_i\,f_j 
+ 2|f_{ij}|^2.\tag 2.9
$$
Therefore we can now conclude that $f_{ij} = 0$.  In particular, 
because
$$
|f_{ij}|^2 \ge |\nabla |\nabla f||^2,
$$
we conclude that $|\n f|$ is identically constant which must be 
identically 0 as $M$ has infinite volume. Therefore,
$f$ is a constant, which is a contradiction.

Assume now that $M$ is an irreducible locally symmetric space 
of rank at least 2 with one nonparabolic end and at least one 
parabolic end. We first observe that Theorem 2.1 is applicable
here. Indeed, for the geodesic distance function $r(x)$ to a fixed
point $p$ on $\tilde M,$ we may choose a parallel orthonormal frame
$\{e_1,e_2,\dots,e_n\}$ along the normal geodesic $\gamma(t)$ 
from $p$ to $x$ such that $e_1=\gamma'(t)$ and $\{e_2,\dots,e_n\}$
diagonizes the curvature tensor $R_{1\alpha 1 \mu}$ with 
corresponding eigenvalues $-b_\alpha^2,$ $\alpha=2,\dots,n.$
Then it is easy to see that

$$
\d r=\sum_{\alpha=2}^n b_\alpha\,\coth (b_\alpha \, r).
$$
Also, one computes that

$$
\lambda_1(\tilde M)
=\frac{\left(\sum_{\alpha=2}^n b_\alpha\right)^2}{4}.
$$
Now it is not difficult to see that Theorem 2.1 can be applied to
$M$ with 

$$
f(r)=\sum_{\alpha=2}^n b_\alpha\, \coth (b_\alpha \, r)
$$
and

$$
2a=\sum_{\alpha=2}^n b_\alpha.
$$
Thus, $M$ has no finite volume end or 

$$
\d \beta=a^2,
$$

$$
|\nabla \beta|=1.
$$

In the latter case, following the argument as in Theorem 1.1, 
we fix a level set $N_0$ of $\beta$ and consider a geodesic 
$\tau$ given by $\tau' = \n \beta$ with $\tau(0) \in N_0.$  
At the point $\tau(0),$ let us consider the curvature
$R_{1\alpha 1 \mu}$ as a bilinear form restrict on the tangent 
space of $N_0$.  In particular, there exists an orthonormal
frame $\{e_1, e_2, \dots e_n\}$ with $e_1= \tau'$ 
and $e_\alpha \in TN_0$ for all $2 \le \alpha \le n$ such that 
$\{e_2,\dots,e_n\}$ diagonizes $R_{1\alpha 1 \mu}$.  
Since $M$ is an irreducible locally symmetric manifold of noncompact 
type, the sectional curvature of $M$ must be nonpositive, hence 

$$
R_{1\alpha 1 \mu}= -b_{\alpha}^2\,\delta_{\alpha \mu}
$$
for some $b_\alpha \ge 0$.  We extend the orthonormal frame 
along $\tau$ by parallel translating the basis 
$\{e_1, e_2, \dots, e_n\}.$  
Since $M$ is locally symmetric, the curvature satisfies
$$
\split
\frac{\p R_{1\alpha 1 \mu}}{\p t} &= R_{1 \alpha 1 \mu, 1}\\&=0.
\endsplit
$$
So
$$
R_{1\alpha 1 \mu} = -b^2_{\alpha} \, \delta_{\alpha \mu}
$$
along $\tau.$  We consider the vector field
$$
V_\alpha(t) = \exp(-b_{\alpha} t)\,e_\alpha
$$
and verify that
$$
\n_{\tau'} \n_{\tau'} V_\alpha = b^2_{\alpha}\,V_\alpha.
$$
On the other hand, 
$$
\split
R_{\tau' V_\alpha}\tau' 
&= \exp(-b_{\alpha}t)\,R_{1 \alpha 1 \alpha} \,e_\alpha\\
&= -b^2_{\alpha}\,V_\alpha.
\endsplit
$$
Hence, $V_\alpha$ is a Jacobi field along $\tau.$  Since this is 
true for all $2 \le \alpha \le n,$ we conclude that the metric 
on $N_t$ must be of the form
$$
ds_t^2 = \sum_{\alpha=2}^n \exp(-2b_\alpha t)\, \omega_\alpha^2,\tag2.10
$$
where $\{ \omega_\alpha\}_{\alpha=2}^n $ is the dual coframe to 
$\{e_\alpha\}_{\alpha=2}^n$ 
at $N_0.$ In particular, the second fundamental form on $N_t$ must 
be a diagonal matrix when written in terms of the 
basis $\{e_\alpha\}_{\alpha=2}^n$ with eigenvalues given by $\{-b_\alpha\}_{\alpha=2}^n$.  
Moreover, since the Buseman function $\beta$ has no critical points, 
and for any $x_0 \in N$, the curve $(t, x_0)$ is a geodesic line 
which can be used to define $\beta,$ the level sets $N_t$ must be 
a compact quotient of some horosphere on $\tilde M.$
\qed
\enddemo

\heading
\S3 Nonparabolic Ends
\endheading

In this section, we prove Theorem B. We begin with the following result.

\proclaim{Theorem 3.1} Let $M^m$ be a complete K\"ahler manifold 
with complex dimension $m \ge 2.$  Assume that $M$ satisfies 
property ($\Cal P_\rho$)
for some nonzero weight function $\rho \ge 0.$  
Suppose the Ricci curvature of $M$ is bounded below by
$$
\text{Ric}_M(x) \ge - 4 \rho(x)
$$
for all $x \in M.$  If $\rho$ satisfies the property

$$
\lim_{x\to \infty} \rho(x) =0,
$$
then either
\roster
\item $M$ has only one nonparabolic end; or
\item $M$ has two nonparabolic ends and it is diffeomorphic 
to $\Bbb R \times \bar M$ for some compact manifold $\bar M.$  
Moreover, the universal covering $\tilde M$ of $M$ is given by 
the total space of a fiber bundle 
$N^{m-1} \rightarrow \tilde M \rightarrow \Sigma$ 
over a Riemann surface $\Sigma$.  Futhermore, the fibers  
are totally geodesic, holomorphic submanifolds given by $N^{m-1}$.
\endroster
\endproclaim

\demo{Proof} Let us assume that $M$ has at least two nonparabolic 
ends and hence by the theory of Li-Tam \cite{L-T}, there exists a 
bounded harmonic function $f$ with finite Dirichlet integral 
constructed as in \cite{L-W1}. We may assume that $\inf f = 0$ 
and $\sup f = 1$ with the infimum of $f$ achieving at infinity 
of a nonparabolic end $E$ and the supremum of $f$ achieving at 
infinity of the other nonparabolic end given by $M \- E.$ 
In fact, since $f$ has finite Dirichlet integral (see \cite{L-T} 
and \cite{L-W1}), $f$ must be pluriharmonic \cite{L} and satisfies 
the improved Bochner formula (see \cite{L-W2})
$$
\d g \ge -2\rho \,g + g^{-1}\,|\n g|^2\tag 3.1
$$
with $g = |\n f|^{\frac{1}{2}}.$

We will first show that inequality \thetag{3.1} is in fact equality.
To see this, let us consider $\phi$ to be a non-negative  
Lipschitz function with compact support in $M$. Then
$$
\int_M |\n (\phi\,g)|^2= \int_M |\n \phi|^2\,g^2 
+ 2 \int_M \phi\,g\,\la
\n \phi, \n g\ra +\int_M \phi^2\,|\n g|^2.\tag 3.2
$$
The second term on the right hand side can be written as
$$
\split
\int_M \phi\,g\,\la \n \phi, \n g\ra
&= \frac 14\int_M \la \n (\phi^2),\n (g^2)\ra\\
& = - \frac 12 \int_M  \phi^2\,g\,\d g 
- \frac 12\int_M \phi^2\,|\n g|^2\\
&=  \int_M \phi^2\,\rho\,g^2- \int_M \phi^2\,|\n g|^2  
-\frac 12\int_M \phi^2\,g\,h,
\endsplit
$$
where
$$
h = \d g+2\rho \,g - g^{-1}\,|\n g|^2.
$$
Combining with \thetag {3.2} we have
$$
\split
\int_M |\n (\phi\,g)|^2 &+ \frac 12\int_M \phi^2\, g\, h \\
&=\int_M \phi^2\,\rho\,g^2 + \int_M |\n \phi|^2\,g^2 
+ \int_M \phi\,g\,\la \n \phi, \n g\ra.
\endsplit
$$
By the weighted Poincar\'e inequality, this becomes
$$
\frac 12 \int_M \phi^2\, g\, h \le  \int_M |\n \phi|^2\,g^2 
+ \int_M \phi\,g\,\la \n \phi, \n g\ra.\tag3.3
$$

Let us choose $\phi = \psi\,\chi$ to be the product of two 
compactly supported functions.  For $0< \delta <1$ and 
$0< \epsilon < \frac 12$, we set $ \chi$ to be
$$
\chi(x) = \left\{ \aligned 0 \qquad \qquad \quad 
& \qquad \text{on} \qquad 
\Cal L(0, \delta \epsilon) \cup \Cal L(1-\delta\epsilon , 1) \\
(-\log \delta)^{-1} (\log f - \log (\delta \epsilon))  
& \qquad \text{on} \qquad \Cal L(\delta\epsilon, \epsilon) \cap E\\
(-\log \delta)^{-1} (\log (1-f) - \log (\delta \epsilon)) 
& \qquad \text{on} 
\qquad \Cal L(1-\epsilon, 1- \delta\epsilon) \cap (M \- E)\\
1 \qquad \qquad \quad &  \qquad \text{otherwise,}
\endaligned \right.
$$
where $\Cal L(a, b)= \{x \in M\,|\, a \le f(x) \le b\}.$
For $R >0$, we set
$$
\psi(x) = \left\{ \aligned 1\quad 
& \qquad \text{on} \qquad B_\rho(R-1)\\
R-r_\rho  & \qquad \text{on} \qquad B_\rho(R) \- B_\rho(R-1)\\
0 \quad & \qquad \text{on} \qquad M\-B_\rho(R), \endaligned \right.
$$
where $r_\rho$ is the geodesic distance function associated to 
the metric $ds^2_\rho.$
Applying to the first term on the right hand side of \thetag{3.3},
we obtain
$$
\int_M |\n \phi|^2 \, g^2 \le 2\int_M |\n \psi|^2\, \chi^2 \, |\n f|
+ 2\int_M |\n \chi|^2 \, \psi^2 \, |\n f|.\tag 3.4
$$

By the assumption on the Ricci curvature, the local gradient 
estimate of Cheng-Yau \cite{C-Y} (see \cite{L-W2})
for positive harmonic functions asserts that for all $R_0>0$,
$$
|\n f|(x) \le \left((2m-1)\sup_{B(x, R_0)}\sqrt{\rho(y)} 
+ C\,R_0^{-1} \right) \,f(x),\tag 3.5
$$
where $C$ is a constant depending only on $n,$ and $B(x, R_0)$ 
is the ball of radius $R_0$ centered at $x$ with respect to the 
background metric $ds_M^2.$
Let us now choose $R_0 = (\sup_{B(x, R_0)} \sqrt{\rho})^{-1}.$  
This choice of $R_0$ is possible as the function 
$r - (\sup_{B(x, r)} \sqrt{\rho})^{-1}$
is negative when $r \to 0$ and it tends to $\infty$ as $r \to \infty.$
Let us observe that if $y \in B(x, R_0)$,  and if $\gamma$ 
is a $ds_M^2$ minimizing geodesic joining $x$ to $y$, then
$$
\split
r_\rho(x, y) &= \int_{\gamma} \sqrt{\rho(\gamma(t)}\, dt\\
&\le R_0\,\sup_{B(x, R_0)} \sqrt{\rho(y)}\\
&\le 1.
\endsplit
$$
This implies that $B(x, R_0) \subset B_\rho(x,1).$ 
Hence \thetag{3.5} can be written as
$$
|\n f|(x) 
\le C\, \left(\sup_{B_\rho(x, 1)} \sqrt{\rho}\right)\, f(x). \tag 3.6
$$
Similarly, applying the same estimate to $1-f$, we also have
$$
|\n f|(x) \le C\, 
\left(\sup_{B_\rho(x, 1)} \sqrt{\rho}\right)\, (1-f(x)).\tag 3.7
$$

At the end $E$, the first term on the right hand side 
of \thetag{3.4} can be estimated by
$$
\split
\int_{E} |\n \psi|^2\, \chi^2 \, |\n f| 
&\le \int_{\Omega} \rho\,|\n f|\\
&\le \left(\int_{\Omega} |\n f|^2 \right)^{\frac{1}{2}} 
\left( \int_{\Omega} \rho^{2}\right)^{\frac 1{2}},
\endsplit \tag3.8
$$
where
$$
\Omega = E \cap (B_\rho(R) \- B_\rho(R-1)) \cap 
\Cal L(\delta \epsilon, 1-\delta\epsilon) .\tag3.9
$$
Recall that \thetag{2.10} and Corollary 2.3 of \cite{L-W4} assert that
$$
\int_{E \cap (B_\rho(R) \- B_\rho(R-1))} \rho\, f^2 
\le C\, \exp(-2R) \tag3.10
$$
and
$$
\int_{E \cap (B_\rho(R) \- B_\rho(R-1))} |\n f|^2  
\le C\, \exp(-2R).\tag3.11
$$
This implies that
$$
\split
\int_{\Omega} \rho^{2}
& \le S^{2}(B_\rho(R)) \,  \int_{\Omega} \rho\\
&   \le S^{2}(B_\rho(R)) \, (\delta\epsilon)^{-2}\, 
\int_{\Omega} \rho\, f^2\\
& \le C\, S^{2}(B_\rho(R)) \,  (\delta\epsilon)^{-2}\, \exp(-2R),
\endsplit
$$
where
$$
S(B_\rho(R)) = \sup_{B_\rho(R)} \sqrt{\rho}.
$$
Hence together with \thetag{3.8} and \thetag{3.11}, we obtain
$$
\int_{E} |\n \psi|^2\, \chi^2 \, |\n f| 
\le C\,S(B_\rho(R)) \, (\delta\epsilon)^{-1}\, \exp(-2R).\tag 3.12
$$

Using \thetag{3.6}, the second term on the right hand side 
of \thetag{3.4} at $E$ can be estimated by
$$
\split
\int_{E}& |\n \chi|^2 \, \psi^2 \, |\n f|\\ & \le (-\log \delta)^{-2}\,
\int_{\Cal L(\delta \epsilon, \epsilon) \cap E \cap B_\rho(R)} 
|\n f|^{3}\,f^{-2}\\
&\le C\, S(B_\rho(R+1))\, (-\log \delta)^{-2}\,
\int_{\Cal L(\delta \epsilon, \epsilon) \cap E \cap B_\rho(R)} 
|\n f|^2\,f^{-1}.
\endsplit\tag3.13
$$
Note that the co-area formula and Lemma 5.1 of \cite{L-W4} imply that
$$
\split
\int_{\Cal L(\delta \epsilon, \epsilon) \cap E \cap B_\rho(R)}
& |\n f|^2\,f^{-1}\\
&\le \int_{\delta \epsilon}^{\epsilon} t^{-1} 
\int_{\ell(t) \cap E_1 \cap B_\rho(R)} |\n f|\,dA\,dt\\
&\le \int_{\ell(b)} |\n f|\,dA 
\int_{\delta \epsilon}^{ \epsilon} t^{-1}\, dt
\endsplit
$$
for a fixed level $b,$ where $\ell(t) = \{ x \in M\, |\, f(x) = t\}.$
Since
$$
\int_{\delta \epsilon}^{ \epsilon} t^{-1}\, dt
 = - \log \delta,
$$
together with \thetag{3.13}, we conclude that there exists 
a constant $C_1 >0,$ such that,
$$
\split
\int_{E} &|\n \chi|^2 \, \psi^2 \, |\n f| \\
& \le
C_1\, S(B_\rho(R+1))\,(-\log \delta)^{-1}.
\endsplit \tag3.14
$$
Combining with \thetag{3.4} and \thetag{3.12}, we get
$$
\int_E |\n \phi|^2 \, g^2 \le C_2\,S(B_\rho(R+1))\,
\left( \delta^{-1} \epsilon^{-1}\, \exp(-2R) 
+ (-\log \delta)^{-1}\right).\tag3.15
$$
By first letting $R \to \infty$ and then letting $\delta \to 0,$ 
the right hand side of \thetag{3.15} vanishes on $E$.  
A similar estimate also works on $M \- E$ by using the 
function $1-f$ instead.

To estimate the second term on the right hand side of \thetag{3.3}, 
we consider
$$
2\int_M \phi\,g\, \la \n \phi, \n g\ra 
= \int_M \psi\, \chi^2 \, \la \n \psi, \n (g^2)\ra 
+ \int_M \psi^2\, \chi\,\la \n \chi, \n (g^2) \ra.\tag 3.16
$$
On $E,$ the first term on the right hand side can be estimated by
$$
\split
\int_E \psi\, \chi^2 \la \n \psi, \n (g^2) \ra 
&\le \int_E \psi \, \chi^2 |\n \psi|\, |\n (g^2)|\\
& \le \int_{\Omega} \sqrt{\rho} \,|\n (g^2)|\\
& \le \left(\int_{\Omega} \rho \right)^{\frac 12} 
\left( \int_{\Omega} |\n (g^2)|^2 \right)^{\frac12},
\endsplit \tag3.17
$$
where
$\Omega$ is given by \thetag{3.9}.  Note that  
\thetag{3.10} asserts that
$$
(\delta\, \epsilon)^2\, \int_{\Omega} \rho \le C\, \exp(-2R).\tag 3.18
$$
Also, since the Bochner formula \thetag{3.1} implies that
$$
\d (g^2) \ge - 4 \rho\, g^2,
$$
if $\tau$ is a nonnegative compactly supported function, then
$$
\split
- 4\int_M \tau^2\, \rho\, g^4 &\le \int_M \tau^2\, (g^2) \d (g^2)\\
& = -2 \int_M \tau\,g^2\, \la \n \tau, \n (g^2) \ra 
- \int_M \tau^2\, |\n (g^2)|^2\\
& \le 2 \int_M |\n \tau|^2\, g^4 
- \frac 12 \int_M \tau^2\, |\n (g^2)|^2,
\endsplit
$$
hence
$$
\int_M \tau^2\, |\n (g^2)|^2 \le 8 \int_M \tau^2\, \rho\, g^4 
+ 4\int_M |\n \tau|^2\,g^4.\tag3.19
$$
Let us set
$$
\tau = \left\{ \aligned 0 \qquad & \qquad \text{on} 
\qquad B_\rho(R-2) \cup (M \- B_\rho(R+1))\\
r_\rho - R +2 & \qquad \text{on} \qquad B_\rho(R-1) \- B_\rho (R-2)\\
1 \qquad & \qquad \text{on} \qquad B_\rho(R) \- B_\rho(R-1)\\
R-r_\rho+1  & \qquad \text{on} \qquad B_\rho(R+1)\- B_\rho(R). 
\endaligned \right.
$$
Then \thetag{3.19} implies
$$
\split
\int_{B_\rho(R) \- B_\rho(R-1)} |\n (g^2)|^2 
&\le C\, \int_{B_\rho(R+1)\- B_\rho(R-2)} \rho\, g^4\\
&\le C\,S^2(B_\rho(R+1)) \int_{B_\rho(R+1)\- B_\rho(R-2)} |\n f|^2.
\endsplit \tag3.20
$$
Applying \thetag{3.11} to \thetag{3.20} and combining with 
\thetag{3.17} and \thetag{3.18}, we conclude that
$$
 \int_E \psi\, \chi^2 \la \n \psi, \n (g^2) \ra 
\le C\, (\delta\, \epsilon)^{-1}\,S(B_\rho(R+1))\, \exp(-2R).\tag3.21
$$

To estimate the second term on the right hand side of \thetag{3.16} on $E,$ 
we integrate by parts and get
$$
\split
\int_E \psi^2\, \chi\, \la \n \chi, \n (g^2) \ra 
& =\int_{\Cal L(\delta \epsilon, \epsilon)\cap E} 
\psi^2\, \chi\, \la \n \chi, \n (g^2) \ra\\
&= -\int_{\Cal L(\delta\epsilon, \epsilon)\cap E} 
\psi^2\, \chi\, \d \chi \,g^2
- \int_{\Cal L(\delta \epsilon, \epsilon)\cap E } 
\psi^2\, g^2\, |\n \chi|^2 \\
& \qquad - 2 \int_{\Cal L(\delta \epsilon, \epsilon)\cap E} 
\psi\, \chi \la \n \psi, \n \chi \ra\, g^2 
+ \int_{\ell(\epsilon)\cap E} \psi^2\, \chi\,\chi_\nu\, g^2 \\
& \qquad - \int_{\ell(\delta \epsilon)\cap E} 
\psi^2\, \chi\, \chi_\nu\, g^2,
\endsplit\tag3.22
$$
where $|\n f|\,\nu = \n f.$
Using the definition of $\chi$ and \thetag{3.6}, 
the two boundary terms can be estimated by
$$
\split
 \int_{\ell(\epsilon)\cap E} \psi^2\, \chi\,\chi_\nu\, g^2 
- \int_{\ell(\delta \epsilon)\cap E } \psi^2\, \chi\, \chi_\nu\, g^2
&= (-\log \delta)^{-1}\,\int_{\ell(\epsilon)\cap E} 
\psi^2 f_\nu\,f^{-1}\,g^2\\
& \le C\,S(B_\rho(R+1))\,(-\log \delta)^{-1}\,
\int_{\ell(\epsilon)} \psi^2 |\n f|\\
& \le C\,S(B_\rho(R+1))\,(-\log \delta)^{-1}.
\endsplit \tag 3.23
$$
We can also estimate  the term
$$
\split
-\int_{\Cal L(\delta\epsilon, \epsilon)\cap E} 
\psi^2\, \chi\,  \d \chi \, g^2
& = (-\log \delta)^{-2} 
\int_{\Cal L(\delta \epsilon, \epsilon) \cap B_\rho(R)\cap E} 
\psi^2\, g^2\,(\log f - \log \delta \epsilon) |\n f|^2\, f^{-2}\\
& \le C\,S(\bar \Omega)\,(-\log \delta)^{-2} 
\int_{\Cal L(\delta \epsilon, \epsilon) \cap B_\rho(R)\cap E} 
\psi^2\,(\log f - \log \delta \epsilon) |\n f|^2\,f^{-1}\\
& \le C\,S(\bar \Omega)\,(-\log \delta)^{-2} 
\int_{\delta \epsilon}^{\epsilon} 
(\log t - \log \delta \epsilon)\,t^{-1} \,dt\\
& \le \frac C2 \, S(\bar \Omega),
\endsplit \tag 3.24
$$
where $S(\bar \Omega) = \sup_{\bar \Omega} \sqrt{\rho}$ 
with $\bar \Omega = \{ x \in M \, |\, r_\rho(x, B_\rho(R) 
\cap \Cal L (\delta \epsilon, \epsilon)) \le 1 \}.$
Finally, Schwarz inequality implies that
$$
-\int_{\Cal L(\delta \epsilon, \epsilon)\cap E} \psi^2\,g^2\,|\n  \chi|^2 
- 2 \int_{\Cal L(\delta \epsilon, \epsilon)\cap E} 
\psi\, \chi\, \la \n \psi, \n \chi\ra\,g^2 
\le \int_{\Cal L (\delta \epsilon, \epsilon)} \chi^2\,g^2 |\n \psi|^2.
$$
Hence combining with\thetag{3.12},  \thetag{3.23} and \thetag{3.24}, 
we deduce that \thetag{3.22} has the estimate
$$
\int_E \psi^2\, \chi\, \la \n \chi, \n (g^2) \ra 
\le C\, S(B_\rho(R+1)) (-\log \delta)^{-1} 
+ \frac C2 \, S(\bar \Omega).\tag 3.25
$$
Note that since $f$ is harmonic with $ \inf_M f = 0,$ the maximum 
principle asserts that for any fixed compact set $\Omega$, 
the set $\Cal L(\delta \epsilon, \epsilon)$ must not 
intersect $\Omega$ when $\epsilon$ is sufficiently small.  
Hence using the assumptions that $ds_\rho^2$ is complete and
$$
\lim_{x \to \infty} \rho(x) = 0,
$$
we conclude that
$$
S(\bar \Omega) \to 0
$$
as $\epsilon \to 0.$ This implies that \thetag{3.25} becomes
$$
\int_E \psi^2\, \chi\, \la \n \chi, \n (g^2) \ra 
\le C\, S(B_\rho(R+1))  (-\log \delta)^{-1}
$$
which tends to $0$ on $E$ by letting $\delta \to 0.$  
Again, a similar argument yields the vanishing of this 
term on $M \-E.$  This proves that $h$ must be identically 0.

Tracing through the proof of \thetag{3.1}, the equality 
of \thetag{3.1} implies that all the inequalities used in 
the proof must be equalities.  In particular, 
if $\{e_1, e_2, \dots , e_{2m}\}$ is an orthonormal frame such that
$$
|\n f|\, e_1 = \n f\tag 3.26
$$
and
$$
J e_1 = e_2,\tag 3.27
$$
then
the Hessian of $f$ must satisfy
$$
f_{\alpha \beta}(x) = 0\tag 3.28
$$
for all $ 3\le \alpha, \beta \le 2m$ and for all $x$ 
such that $\n f(x) \neq 0.$  We also have
$$
f_{11} = - f_{22}.
$$
Moreover, $\text{Ric}_M (e_1, e_1)(x) = -\rho(x)$ for all $x \in M.$

Since $f$ is pluriharmonic, locally $f$ can be taken to be the 
real part of holomorphic function $F.$  In fact, when lifted to 
the universal covering $\tilde M$ of $M$, $F$ is globally defined.  
The Cauchy-Riemann equations, \thetag{3.26}, and \thetag{3.27}  
imply that the level set of $F$ is orthogonal to $e_1$ and $e_2.$  
Moreover, \thetag{3.28} asserts that the second fundamental form 
of the level set of $F$ with respect to $e_1$ is identically 0.  
Taking $J$ of this, we conclude that the second fundamental form 
of the level set of $F$ with respect to $e_2$ is also identically 0.  
Hence, the level sets of $F$ are totally geodesic at noncritical 
points of $f.$  This gives a totally geodesic holomorphic fibration 
of $\tilde M$ with fibers given by the level sets of $F$ denoted 
by $N.$
\qed
\enddemo

The following theorem deals with the case when $\rho$ is 
just bounded.  Without the help of $\rho \to 0$ at infinity, 
we need to impose a 
stricter assumption on the curvature.

\proclaim{Theorem 3.2} Let $M^m$ be a complete K\"ahler manifold 
with complex dimension $m \ge 2.$  Assume that $M$ satisfies 
property ($\Cal P_\rho$)
for some weight function $\rho \ge 0.$  Suppose the Ricci 
curvature of $M$ is bounded below by
$$
\text{Ric}_M(x) \ge - 4(1-\epsilon) \rho(x)
$$
for some $\epsilon >0$ and for all $x \in M.$  If $\rho$ satisfies

$$
\liminf_{R\to \infty} \frac{S(B_\rho(R))}{R}=0,
$$
where $S(B_\rho(R))=\sup_{B_\rho(R)}\sqrt \rho,$ then $M$ 
must have only
one nonparabolic end.
\endproclaim

\demo{Proof} Following a similar argument as in the proof of 
Theorem 3.1 we assume that $M$ has at least two nonparabolic 
ends and we construct a bounded harmonic function $f.$  Again 
using $g = |\n f|^{\frac 12},$ we have a Bochner formula similar 
to \thetag{3.1}, which  now takes the form
$$
\d g \ge -2(1-\epsilon)\,\rho\,g + g^{-1}\, |\n g|^2.
$$
Inequality \thetag{3.3} will now take the form
$$
\int_M \phi\,g\,h  + 2\epsilon \int_M \phi^2\, \rho\,g^2 \le
2\int_M |\n \phi|^2\,g^2 + 2\int_M \phi\,g\, \la \n \phi, \,\n
g\ra,
$$
with $h = \d g + 2(1-\epsilon)\rho\, g - g^{-1}|\n g|^2$ which is
nonnegative. Hence we conclude that
$$
\split
\epsilon \int_M \phi^2\, \rho\,g^2 &\le \int_M |\n \phi|^2\,g^2 
+ \int_M \phi\,g\, \la \n \phi, \, g \ra\\
&\le (1+\frac 2\epsilon)\int_M |\n \phi|^2\,g^2 
+ \frac \epsilon 2 \int_M \phi^2\,\rho\,g^2
\endsplit
$$
or
$$
\int_M \phi^2\,\rho\,g^2 
\le \frac 2\epsilon (1+\frac 2\epsilon)\int_M |\n \phi|^2\,g^2.
$$
To get a contradiction, we need only to show that the right hand 
side tends to 0 with the appropriate cut-off function $\phi.$  
Using the same cut-off function $\phi = \chi\, \psi$ as before 
and applying the estimates \thetag{3.12} and \thetag{3.14}, we 
conclude that
$$
\int_M \n \phi|^2\,g^2 
\le C\,S(B_\rho(R+1))\,((\delta \epsilon)^{-1}\,
\exp(-2R) + (-\log \delta)^{-1}) \tag 3.29
$$
where $C$ is a constant independent of $\rho.$  First we choose 
a sequence 
$R_i\to \infty$ with the property that
$$
q(R_i)=\sqrt{\frac{S(B_\rho(R_i+1))}{R_i}}\to 0
$$
and then by setting $\delta=\epsilon=\exp(-R_i\,q(R_i)),$ the right
hand side of \thetag {3.29} goes to $0.$ The theorem now follows. 
\qed
\enddemo

Theorem B now follows as a corollary  of Theorem 3.2.  

\proclaim{Corollary 3.3} Let $M^m$ be a complete K\"ahler manifold 
satsifying
$$
\text{Ric}_M > -\frac{\lambda_1(M)}{4}.
$$
Then $M$ must have only one infinite volume end.
\endproclaim

We now claim that the value $\frac{\lambda_1(M)}4$ is best possible 
due to the following class of examples.
Consider manifolds of the form $M = N \times \Sigma,$ where
$N^{m-1}$ is a compact $(m-1)$-dimensional K\"ahler manifold with
Ricci curvature bounded from below by $\text{Ric}_N \ge -1$ and
$\Sigma$ a Riemann surface with constant sectional curvature
$-1.$  Moreover, we assume that $\lambda_1(\Sigma) =\frac 14$. It
is easy to see that $\lambda_1(M) = \frac 14$ and $\text{Ric}_M
\ge -1.$   However, it is known \cite{L-W2} that $\Sigma$ can have
more than one infinte volume ends.  One such case is $\Sigma = \Bbb R
\times \Bbb S^1$ with the warped product metric given by
$ds_{\Sigma}^2 = dt^2 + \coth^2 t \, d\theta^2.$  This example
shows that the bound on the Ricci curvature in Theorem 3.2 and 
Corollary 3.3 is sharp.

\heading
\S4 Parabolic Ends with $\rho \to 0$
\endheading

In this section, we will consider the corresponding theorem to 
Theorem A for K\"ahler manifolds with property ($\Cal P_\rho$).

\proclaim{Theorem 4.1} Let $M^m$ be a complete K\"ahler manifold 
of dimension $m \ge 2$ with property ($\Cal P_\rho$).  Suppose 
the holomorphic bisectional curvature of $M$ is bounded from below by
$$
R_{i \bar i j \bar j} (x) \ge -\frac{\rho(x)}{m^2} \qquad 
\text{for all} \qquad i \neq j,
$$
$x \in M$, and for all unitary frame $\{e_1, e_2, \dots e_m\}.$  
If $\rho$ satisfies the property that
$$
\lim_{x\to \infty} \rho(x)=0,
$$
then $M$ has at most 2 ends if $m \ge 3$.

In the case when $m=2$, $M$ has at most 4 ends.  Moreover, 
if $M$ has exactly 4 ends, then
$$
R_{i \bar i j \bar j} (x) = -\frac{\rho(x)}{4} \qquad 
\text{for all} \qquad i\neq j.
$$
\endproclaim

\demo{Proof} Note that by \cite{L-W4} the existence of the weight function $\rho$ 
asserts that $M$ must be nonparabolic.  On the other hand, 
Theorem 3.2 implies that $M$ must have only one nonparabolic end.  
If $M$ has $k$ ends, with $k >1,$ then $M$ must have $k-1$ parabolic 
ends.  The theory of Li-Tam \cite{L-T} asserts that for each 
parabolic end $E$, one can construct a positive harmonic 
function $f$ such that
$$
\limsup_{x\to \infty}f(x)=\infty \qquad \text{with} \qquad x \in E  
$$
and $f$ is bounded on $M \- E$ and
$$
\liminf_{x\to \infty} f(x) = 0 
\qquad \text{with} \qquad x \in M \- E.
$$

According to Lemma 4.1 of \cite{L-W3}, if we let 
$g =|f_{\alpha \bar \beta}|^{\frac 12}$ then it satisfies the 
inequality
$$
\d g \ge -\frac 2m\,\rho\, g - \frac{m-2}{m}\,g^{-1}|\n g|^2.\tag 4.1
$$
Similar to Theorem 3.1, we will show that this inequality is indeed 
an equality.
Following the argument of Theorem 3.1, by applying the weighted 
Poincar\'e inequality, we obtain
$$
\frac m2 \int_M \phi^2\, g\, h 
\le -(m-2) \int_M \phi\, g\,\la \n \phi, \n g \ra 
+ \int_M |\n \phi|^2\, g^2, \tag 4.2
$$
where
$$
h = \d g + \frac 2m\, \rho \, g + \frac{m-2}{m} \, g^{-1}|\n g|^2.
$$
We need to choose the cut-off function $\phi$ so that the right hand 
side of \thetag{4.2} tends to 0 and conclude that $h$ must be 
identically 0.

Note that by a theorem of Nakai \cite{N} (also see \cite{N-R}), 
the positive harmonic function $f$ can be taken to be proper on 
the parabolic end $E$, i.e.,
$$
\liminf_{x \to \infty} f(x) = \infty
$$
for $x \in E.$
To deal with the right hand side of \thetag{4.2} on $E$,   we define
$$
\phi = \left\{ \aligned 1 \qquad \qquad & 
\qquad \text{on} \qquad \Cal L (0, 2T) \cap E\\
T^{-1}(3T-f)& \qquad \text{on} \qquad \Cal L (2T, 3T) \cap E\\
0 \qquad \qquad 
& \qquad \text{on} \qquad \Cal L(3T, \infty) \cap E.
\endaligned \right.
$$
After integrating by parts, we have
$$
\split
-\int_E \phi\,g\, \la \n \phi\, \n g \ra 
&= -\frac 14 \int_{\Cal L(2T, 3T) \cap E}  \la \n \phi^2, \n g^2 \ra\\
& = \frac 14 \int_{\Cal L(2T, 3T) \cap E} \d (\phi^2) g^2   
+ \frac 12 \int_{\ell (2T) \cap E} \phi_\nu\, g^2,
\endsplit\tag4.3
$$
where $\nu$ is the unit normal to the level set $\ell(2T)$ given 
by $|\n f|\, \nu = \n f.$  Using the definition of $\phi$, we obtain
$$
\split
\int_{\Cal L(2T, 3T) \cap E} \d (\phi^2)\, g^2
&= 2\int_{\Cal L(2T, 3T) \cap E} |\n \phi|^2\, g^2 
+ 2 \int_{\Cal L(2T, 3T)\cap E} \phi\, \d \phi\, g^2\\
& = 2T^{-2} \,\int_{\Cal L(2T, 3T) \cap E} |\n f|^2\, g^2\\
& \le 2T^{-2} \left( \int_{\Cal L(2T, 3T) \cap E} 
|\n f|^4 \right)^{\frac 12}\,
\left(\int_{\Cal L(2T, 3T) \cap E} 
|f_{\alpha \bar \beta}|^2 \right)^{\frac 12}.
\endsplit \tag 4.4
$$

Applying the assumption on the curvature and the gradient 
estimate \thetag{3.6}, we conclude that there exists a 
constant $C >0$ such  that
$$
|\n f|(x) \le  C\, S(N_1(\Cal L(2T, 3T) \cap E))\, f(x)
$$
for all $x\in \Cal L(2T, 3T) \cap E,$ where
$$
S(N_1(\Cal L(2T, 3T) \cap E))
= \sup_{N_1(\Cal L(2T, 3T) \cap E))} \sqrt{\rho}
$$
with the supremum taken over the set
$$
N_1(\Cal L(2T, 3T) \cap E)
= \{ x \, |\, r_\rho(x, \Cal L(2T, 3T) \cap E) \le 1\}.
$$
Also, since
$$
\split
0 &=
\int_{\Cal L(1, t) \cap E} \d f \\
&= -\int_{\ell(1) \cap E} f_{\nu} + \int_{\ell(t) \cap E} f_\nu\\
& = -\int_{\ell(1) \cap E} |\n f| + \int_{\ell(t)\cap E} |\n f|
\endsplit
$$
where $\nu$ is the outward pointing unit normal vector 
to $\ell(1)$ and $\ell(t),$ together with the co-area formula, 
the first integral on the right hand side of \thetag{4.4} can 
be estimated by
$$
\split
\int_{\Cal L(2T, 3T) \cap E} |\n f|^4 
&\le \int_{2T}^{3T} \int_{\ell(t) \cap E} |\n f|^3 \,dA(t)\, dt\\
&\le C\,S^2(N_1(\Cal L(2T, 3T) \cap E)) 
\int_{2T}^{3T} t^2 \int_{\ell (t) \cap E} |\n f|\,dA(t)\, dt\\
& \le C_1\, S^2(N_1(\Cal L(2T, 3T) \cap E))\,T^3.
\endsplit \tag 4.5
$$

We will now estimate the second integral on the right hand side 
of \thetag{4.4}.  Let us observe that the Bochner formula asserts that
$$
\d |\n f|^2 \ge 2|f_{ij}|^2 - 2m^{-2} \rho\,|\n f|^2,
$$
where $(f_{ij})$ is the real Hessian of $f.$
Multiplying with a cut-off function $\varphi^2$ and integrating by 
parts yield
$$
\int_M \varphi^2\, \d |\n f|^2 \ge 2 \int_M \varphi^2\, |f_{ij}|^2 
- 2m^{-2} \int_M \varphi^2\, \rho\, |\n f|^2.\tag 4.6
$$
On the other hand,
$$
\split
\int_M \varphi^2\, \d |\n f|^2 
&= - 2\int_M \varphi\,|\n f|\, \la \n \varphi, \n |\n f| \ra\\
& \le \int_M \varphi^2\, |\n |\n f||^2 
+ \int_M |\n \varphi|^2\, |\n f|^2.
\endsplit
$$
Using the inequality
$$
|f_{ij}|^2 \ge |\n |\n f||^2
$$
and combining with \thetag{4.6}, we obtain
$$
2m^{-2} \int_M \varphi^2\, \rho\, |\n f|^2 
+ \int_M |\n \varphi|^2\, |\n f|^2 
\ge \int_M \varphi^2\, |f_{ij}|^2.\tag4.7
$$
Choosing
$$
\varphi = \left\{ \aligned 0\quad \quad
& \qquad \text {on} \qquad 
\Cal L(0, T) \cup \Cal L(4T, \infty)\cup (M\- E)\\
T^{-1}(f -T)& \qquad \text {on} \qquad \Cal L(T, 2T) \cap E\\
1 \quad \quad & \qquad \text {on} \qquad \Cal L(2T, 3T) \cap E\\
T^{-1}(4T -f)& \qquad \text {on} \qquad \Cal L(3T, 4T) \cap E,
\endaligned \right.
$$
we conclude that
$$
2m^{-2} \int_{\Cal L(T, 4T) \cap E} \rho\, |\n f|^2 
+ T^{-2}\int_{\Cal L(T, 4T)\cap E} |\n f|^4 
\ge \int_{\Cal L(2T, 3T)\cap E} |f_{ij}|^2.\tag4.8
$$
By applying \thetag{4.5} to the second term on the left hand side 
and using
$$
\split
\int_{\Cal L(T, 4T) \cap E} \rho\, |\n f|^2
&\le S^2(\Cal L(T, 4T) \cap E) \, \int_{\Cal L(T, 4T)\cap E} |\n f|^2\\
& \le S^2(\Cal L(T, 4T) \cap E)\, 
\int_T^{4T} \int_{\ell(t)} |\n f|\,dA(t)\, dt\\
&= C_1\,T\,S^2(\Cal L(T, 4T) \cap E),
\endsplit
$$
\thetag{4.8} becomes
$$
C_1\, S^2(N_1(\Cal L(T, 4T) \cap E))\,T
\ge \int_{\Cal L(2T, 3T) \cap E} |f_{ij}|^2.
$$
Combining with \thetag{4.5} and \thetag{4.4}, we obtain
$$
\int_{\Cal L(2T, 3T)\cap E} \d(\phi^2) g^2 
\le C_2 \, S^2(N_1(\Cal L(T, 4T) \cap E)).\tag 4.9
$$

The second term of the right hand side of \thetag{4.3} can be 
estimated by writing it as
$$
\int_{\ell (2T) \cap E} \phi_{\nu}\,g^2 = T^{-1} \,
\int_{\ell(2T)\cap E} |\n f|\,|f_{\alpha \bar \beta}|.
$$
However, the mean value theorem asserts that there exists 
a $T' \in [2T, 3T]$ such that
$$
T^{-1}\,\int_{\Cal L(2T, 3T)\cap E } |\n f|^2\, 
|f_{\alpha \bar \beta}| =  \int_{\ell(T')\cap E} |\n f|\, 
|f_{\alpha \bar \beta}|.
$$
The estimate we use for \thetag{4.4} implies that
$$
T^{-2} \int_{\Cal L(2T, 3T)\cap E } |\n f|^2\, 
|f_{\alpha \bar \beta}| 
\le C_2 \, S^2(N_1(\Cal L(T, 4T) \cap E))\tag 4.10
$$
which tends to 0 as $T \to \infty$ by the assumption on $\rho.$  
Hence we conclude that there exists a sequence of $T' \to \infty$ 
such that
$$
\int_{\ell (2T') \cap E} \phi_{\nu}\,g^2 \to 0.
$$
Together with \thetag{4.9}, we conclude that the right hand 
side of \thetag{4.3} tends to 0 as $T \to \infty.$

The second term on the right hand side of \thetag{4.2} is given by
$$
\int_E |\n \phi|^2\, g^2 = T^{-2} \int_{\Cal L(2T, 3T) \cap E} 
|\n f|^2\, |f_{\alpha \bar \beta}|
$$
and tends to 0 by \thetag{4.10}.

We will now consider \thetag{4.2} on the remaining manifold $M \-E$. 
Note that the
properness of $f$ on $E$ implies that the sublevel set
$\Cal L(0 , 4\epsilon) \cap E = \emptyset$ for sufficiently 
small $\epsilon.$
By taking $M \- E$ as a nonparabolic end, we choose 
$\phi= \psi\, \chi$ as
in Theorem 5.2 of \cite{L-W4}.  In particular, we set
$$
\chi(x) = \left\{ \aligned 0 \quad \quad 
& \qquad \text{on} \qquad \Cal L(0, 2\epsilon) \\
\epsilon^{-1} (f - 2\epsilon)  
& \qquad \text{on} \qquad \Cal L(2\epsilon, 3\epsilon) \\
1 \quad  \quad 
&  \qquad \text{on} \qquad \Cal L(3\epsilon, \infty) \cap (M \- E),
\endaligned \right.
$$
and
$$
\psi(x) = \left\{ \aligned 1\quad 
& \qquad \text{on} \qquad B_\rho(R-1) \cap (M\-E)\\
R-r_\rho  & \qquad \text{on} 
\qquad (B_\rho(R) \- B_\rho(R-1)) \cap (M\-E)\\
0 \quad & \qquad \text{on} \qquad M\-(B_\rho(R) \cup E)
\endaligned \right.
$$
To estimate the first term on the right hand side of \thetag{4.2} on $M\setminus E,$ 
we write
$$
-2\int_{M\setminus E} \phi\,g\, \la \n \phi\, \n g \ra 
= - \int_{M\setminus E} \psi\, \chi^2 \la \n \psi, \n (g^2) \ra 
- \int_{M\setminus E} \psi^2\, \chi \la \n \chi, \n (g^2)\ra, \tag4.11
$$
with
$$
\split
- \int_{M\setminus E} \psi\, \chi^2 \la \n \psi, \n (g^2) \ra 
&\le \int_{M\setminus E} \psi \, \chi^2 |\n \psi|\, |\n (g^2)|\\
& \le \int_{\Omega} \sqrt{\rho} \,|\n (g^2)|\\
& \le \left(\int_{\Omega} \rho \right)^{\frac 12} 
\left( \int_{\Omega} |\n (g^2)|^2 \right)^{\frac12},
\endsplit \tag4.12
$$
where
$$
\Omega = (B_\rho(R) \- B_\rho(R-1)) 
\cap \Cal L(2\epsilon, \infty) \cap (M\- E).
$$
Note that the estimate \thetag{2.9} of \cite{L-W4} asserts that
$$
4\epsilon^2\, \int_{\Omega} \rho \le C\, \exp(-2R).\tag 4.13
$$
Also, since the Bochner formula \thetag{4.1} implies that
$$
\d (g^2) \ge - \frac{4}{m} \rho\, g^2 +\frac 4m |\n g|^2,
$$
if $\tau$ is a nonnegative compactly supported function, then
$$
\split
- 4m^{-1} \int_M \tau^2\, \rho\, g^4 
+ m^{-1} \int_M \tau^2\,|\n (g^2)|^2 
&\le \int_M \tau^2\, (g^2) \d (g^2)\\
& = -2 \int_M \tau\,g^2\, \la \n \tau, \n (g^2) \ra 
- \int_M \tau^2\, |\n (g^2)|^2\\
& \le  \int_M |\n \tau|^2\, g^4.
\endsplit
$$
Hence
$$
\int_M \tau^2\, |\n (g^2)|^2 \le 4 \int_M \tau^2\, \rho\, g^4 
+ m\int_M |\n \tau|^2\,g^4.\tag4.14
$$
Let us set
$$
\tau = \left\{ \aligned 0 \qquad & \qquad \text{on} 
\qquad B_\rho(R-2) \cup (M \- B_\rho(R+1))\\
r_\rho - R +2 & \qquad \text{on} \qquad B_\rho(R-1) \- B_\rho (R-2)\\
1 \qquad & \qquad \text{on} \qquad B_\rho(R) \- B_\rho(R-1)\\
R-r_\rho+1 & \qquad \text{on} \qquad B_\rho(R+1)\- B_\rho(R).
\endaligned \right.
$$
Then \thetag{4.14} implies
$$
\split
\int_{B_\rho(R) \- B_\rho(R-1)} |\n (g^2)|^2 
&\le C\, \int_{B_\rho(R+1)\- B_\rho(R-2)} \rho\, g^4\\
&\le C\,S^2(B_\rho(R+1) \-B_\rho(R-2)) 
\int_{B_\rho(R+1)\- B_\rho(R-2)} |\n f|^2.
\endsplit \tag4.15
$$

Applying Corollary 2.3 of \cite{L-W4} to \thetag{4.15} 
and combining with \thetag{4.12} and \thetag{4.13}, we conclude that
$$
- \int_M \psi\, \chi^2 \la \n \psi, \n (g^2) \ra 
\le C\, \epsilon^{-1}\,S(B_\rho(R+1) \- B_\rho(R-2))\, 
\exp(-2R).\tag4.16
$$

To estimate the second term on the right hand side of \thetag{4.11}, 
we integrate by parts and get
$$
\split
-\int_{M\setminus E} \psi^2\, \chi\, \la \n \chi, \n (g^2) \ra 
& =- \int_{\Cal L(2\epsilon, 3\epsilon)\cap (M\setminus E)} 
\psi^2\, \chi\, \la \n \chi, \n (g^2) \ra\\
&= \int_{\Cal L(2\epsilon, 3\epsilon)\cap (M\setminus E)} \psi^2\, \chi\, \d \chi \,g^2
+ \int_{\Cal L(2\epsilon, 3\epsilon)\cap (M\setminus E) } \psi^2\, g^2\, |\n \chi|^2 \\
& \qquad + 2 \int_{\Cal L(2\epsilon, 3\epsilon)\cap (M\setminus E)} 
\psi\, \chi \la \n \psi, \n \chi \ra\, g^2 
- \int_{\ell(3\epsilon)\cap (M\setminus E)} \psi^2\, \chi\,\chi_\nu\, g^2 \\
& \qquad + \int_{\ell(2\epsilon)\cap (M\setminus E)} \psi^2\, \chi\, \chi_\nu\, g^2,
\endsplit\tag4.17
$$
where $|\n f|\,\nu = \n f.$
Using the definition of $\chi$, the two boundary terms become
$$
\split
- \int_{\ell(3\epsilon)\cap (M\setminus E)} \psi^2\, \chi\,\chi_\nu\, g^2 
+ \int_{\ell(2\epsilon)\cap (M\setminus E)} \psi^2\, \chi\, \chi_\nu\, g^2
&= -\epsilon^{-1}\,\int_{\ell(3\epsilon)\cap (M\setminus E)} \psi^2 f_\nu\,f^{-1}\,g^2\\
& \le 0.
\endsplit
$$
Hence \thetag{4.17} becomes
$$
-\int_{M\setminus E} \psi^2\, \chi\, \la \n \chi, \n (g^2) \ra  
\le 2 \int_{\Cal L(2\epsilon,3\epsilon)} \psi^2\, g^2\, |\n \chi|^2
+  \int_{\Cal L(2\epsilon, 3\epsilon)} 
\chi^2\, |\n \psi|^2\, g^2, \tag4.18
$$
where we have used the fact that
$$
\d \chi = 0 \qquad \text{on} \qquad \Cal L(2\epsilon, 3\epsilon).
$$

We will now estimate the first term on the right hand side 
of \thetag{4.18}.  Using the definition of $\chi$, $\psi$ 
and $g,$ we have
$$
\split
\int_{\Cal L(2\epsilon,3\epsilon)} \psi^2\, g^2\, |\n \chi|^2
&\le \epsilon^{-2} \int_{\Cal L(2\epsilon, 3\epsilon)\cap B_\rho(R)} 
|\n f|^2\, |f_{\alpha \bar\beta}| \\
&\le \epsilon^{-2} 
\left( \int_{\Cal L(2\epsilon, 3\epsilon) \cap B_\rho(R)} 
|\n f|^4 \right)^{\frac 12} 
\left( \int_{\Cal L(2\epsilon, 3\epsilon) \cap B_\rho(R)} 
|f_{\alpha \bar \beta}|^2 \right)^{\frac12}.
\endsplit \tag 4.19
$$
By the gradient estimate, the term
$$
\split
\int_{\Cal L(2\epsilon, 3\epsilon) \cap B_\rho(R)} |\n f|^4 
&\le C\, S^2(N_1(\Cal L(2\epsilon, 3\epsilon) \cap B_\rho(R)))\, 
\int_{\Cal L(2\epsilon, 3\epsilon) \cap B_\rho(R)} f^2\,|\n f|^2\\
&\le C_1\,\epsilon^3\,
S^2(N_1(\Cal L(2\epsilon, 3\epsilon) \cap B_\rho(R)))
\endsplit
$$
and \thetag{4.19} becomes
$$
\int_{\Cal L(2\epsilon,3\epsilon)} \psi^2\, g^2\, |\n \chi|^2
\le C_1\,\epsilon^{-\frac 12} 
S^2(N_1(\Cal L(2\epsilon, 3\epsilon) \cap B_\rho(R))) 
\left( \int_{\Cal L(2\epsilon, 3\epsilon) \cap B_\rho(R)} 
|f_{\alpha \bar \beta}|^2 \right)^{\frac12}.\tag4.20
$$
The second term on the right hand side of \thetag{4.18} can be 
estimated by
$$
\split
\int_{\Cal L(2\epsilon, 3\epsilon)} \chi^2\, |\n \psi|^2\, g^2 
& \le  
\int_{\Cal L(2\epsilon, 3\epsilon) \cap (B_\rho(R) \- B_\rho(R-1))} 
\rho \, |f_{\alpha \bar \beta}|\\
& \le 
\left(\int_{\Cal L(2\epsilon, 3\epsilon)\cap (B_\rho(R)\-B_\rho(R-1))}
\rho^2 \right)^{\frac 12} 
\left( \int_{\Cal L(2\epsilon, 3\epsilon)\cap B_\rho(R)} 
|f_{\alpha \bar \beta}|^2 \right)^{\frac 12}\\
& \le S(B_\rho(R)\-B_\rho(R-1))\,(2\epsilon)^{-1}\,
\exp(-R)\, \left( \int_{\Cal L(2\epsilon, 3\epsilon)\cap B_\rho(R)} 
|f_{\alpha \bar \beta}|^2 \right)^{\frac 12}.
\endsplit\tag4.21
$$

We will now estimate the term involving the complex Hessian.  
Using \thetag{4.7} by setting $\varphi = \tau\,\eta$ with
$$
\tau = \left\{ \aligned 0 \quad \quad 
& \qquad \text{on} \qquad \Cal L(0, \epsilon) 
\cup (\Cal L(4\epsilon, \infty) \cap (M \- E)) \\
\epsilon^{-1} (f - \epsilon)  
& \qquad \text{on} \qquad \Cal L(\epsilon, 2\epsilon) \\
1 \quad  \quad 
&  \qquad \text{on} \qquad \Cal L(2\epsilon, 3\epsilon) \\
\epsilon^{-1}(4\epsilon - f) 
& \qquad \text{on} \qquad \Cal L(3\epsilon, 4\epsilon),
\endaligned \right.
$$
and
$$
\eta(x) = \left\{ \aligned 1\quad \quad
& \qquad \text{on} \qquad B_\rho(R) \cap (M\- E)\\
R+1-r_\rho  
& \qquad \text{on} \qquad (B_\rho(R+1) \- B_\rho(R)) \cap (M \-E)\\
0 \quad \quad
& \qquad \text{on} \qquad M\-B_\rho(R+1), 
\endaligned \right.
$$
we have
$$
\split
\int_{\Cal L(2\epsilon, 3\epsilon) \cap B_\rho(R)} 
|f_{\alpha \bar \beta}|^2 &\le \int_M \tau^2\, \eta^2\,
|f_{\alpha \bar \beta}|^2\\
&\le 2m^{-2} \int_M \tau^2\, \eta^2\,\rho\,|\n f|^2 
+ 2\int_M |\n \tau|^2\, \eta^2\,|\n f|^2 
+ 2\int_M \tau^2\, |\n \eta|^2\, |\n f|^2\\
&\le 2(m^{-2}+1) \int_{\Cal L(\epsilon, 4\epsilon) \cap B_\rho(R+1)} 
\rho\, |\n f|^2 
+ 2\epsilon^{-2} \int_{\Cal L(\epsilon, 4\epsilon) \cap B_\rho(R+1)} 
|\n f|^4 \\
& \le C_1\, S^2(\Cal L(\epsilon, 4\epsilon) \cap B_\rho(R+1))\,
\int_{\Cal L(\epsilon, 4\epsilon) \cap B_\rho(R+1)} |\n f|^2\\
&\qquad + C\, S(N_1(\Cal L(\epsilon, 4\epsilon) \cap B_\rho(R+1)))\,
\epsilon^{-2} \int_{\Cal L(\epsilon, 4\epsilon) \cap B_\rho(R+1)} 
f^2\, |\n f|^2\\
&\le C_2\,\epsilon\, 
S^2(N_1(\Cal L(\epsilon, 4\epsilon) \cap B_\rho(R+1))).
\endsplit
$$
Hence, combining with \thetag{4.18}, \thetag{4.20} 
and \thetag{4.21}, we conclude that
$$
\split
-\int_{M\setminus E} \psi^2\, \chi\, \la \n \chi, \n (g^2) \ra  
&\le C_3\,S^2(N_1(\Cal L(2\epsilon, 3\epsilon) \cap B_\rho(R))) \\
& + C\, \epsilon^{-\frac 12}\,e^{-R}\,
S(N_1(\Cal L(\epsilon, 4\epsilon) \cap B_\rho(R+1)))\,
S(B_\rho(R)\-B_\rho(R-1)).
\endsplit
$$
Together with \thetag{4.16} \thetag{4.11}, the first term on 
the right hand side of \thetag{4.2} can be estimated by
$$
\split
-2\int_{M\setminus E} &\phi\,g\,\la \n \phi, \n g\ra \\
&\le  C\, \epsilon^{-1}\,S(B_\rho(R+1)\-B_\rho(R-2))\, \exp(-2R)
+ C_3\,S^2(N_1(\Cal L(2\epsilon, 3\epsilon) \cap B_\rho(R))) \\
&\qquad+ C_2 \epsilon^{-\frac 12}\,\exp(-R)\, 
S(N_1(\Cal L(\epsilon, 4\epsilon) \cap B_\rho(R+1)))\,
S(B_\rho(R) \- B_\rho(R-1)).
\endsplit
$$
By first letting $R \to \infty$ and then $\epsilon \to 0$, 
the decay property of $\rho$ implies that the right hand side 
tends to 0.

The second term on the right hand side of \thetag{4.2} can be 
estimated by
$$
\int_{M\setminus E} |\n \phi|^2\, g^2 \le 2\int_{M\setminus E} \chi^2\, |\n \psi|^2\,g^2 
+ 2\int_{M\setminus E} |\n \chi|^2\, \psi^2\, g^2,
$$
and these two terms on the right hand side can be estimated the 
same way as is \thetag{4.18}.  Hence this term will also tend to 0.  
In particular, we conclude that equality holds on \thetag{4.1}.  
Moreover, the weight function $\rho$ is given by the holomorphic 
bisectional curvature
$$
-\frac{2}{m}\rho \,|f_{\alpha \bar \beta}|^2
= R_{i\bar i j \bar j}(\lambda_i - \lambda_j)^2,
$$
where $\lambda_\alpha$ are the eigenvalues of 
$(f_{\alpha \bar \beta})$ and
hence must be smooth.

Note that equality of \thetag{4.1} can be written as
$$
\d |f_{\alpha \bar \beta}|^{\frac{m-1}{m}} 
= - \frac{4(m-1)}{m^2} \rho \,
|f_{\alpha \bar \beta}|^{\frac{m-1}{m}}.\tag4.22
$$
Since $|f_{\alpha \bar \beta}|^{\frac{m-1}{m}}$ is nonnegative, 
regularity of the differential equation asserts that it must in fact 
be positive.  Moreover, equality of \thetag{4.1} implies that  
inequalities \thetag{4.6}, \thetag{4.7}, and \thetag{4.8} in 
\cite{L-W4} are all equalities.  Equality for \thetag{4.6} in 
\cite{L-W4} asserts that
for each $\theta$, there exists a constant $a_\theta \in \Bbb C$ 
such that
$$
a_\theta \lambda_{\alpha} 
= \p_{\theta} f_{\alpha \bar \alpha} \tag 4.23
$$
for all $\alpha.$  On the other hand, equality 
for \thetag{4.8} in \cite{L-W4} implies that
for each $\alpha$
$$
\p_{\alpha} f_{\theta \bar \theta} 
= \p_{\alpha} f_{\eta \bar \eta}\tag 4.24
$$
for all $\theta, \,\eta \neq \alpha.$
In particular, when $m \ge 3,$ \thetag {4.23} 
and \thetag {4.24} imply that
$$
a_\alpha \lambda_{\theta} = a_\alpha \lambda_{\eta}
$$
for $\theta, \, \eta \neq \alpha.$  If $a_\alpha \neq 0$, this 
implies that $\lambda_{\theta} = \lambda_{\eta}$ for all 
$\theta, \,\eta \neq \alpha.$  
Using $\sum_{\beta} \lambda_\beta = 0$, we conclude that
$$
\lambda_\alpha = -(m-1) \lambda_\theta.\tag 4.25
$$
If there is another $\theta \neq \alpha$ such 
that $a_\theta \neq 0$ , then the same argument implies that
$$
\lambda_\theta = -(m-1) \lambda_\alpha
$$
contradicting \thetag {4.25} unless  $\lambda_\beta = 0$ for 
all $\beta.$  This implies that $g=0$ at that point, hence it 
must be identically 0, implying $f$ is pluriharmonic and $M$ has 
only one end.

At a fixed point $p\in M$, the only situations left are when
either all $a_\alpha = 0$ or there exists exactly one 
$a_\alpha \neq 0.$  When $a_\alpha = 0$ for all $\alpha,$ 
then \thetag {4.23} implies that 
$\p_{\alpha} f_{\theta \bar \theta} = 0$ for all 
$\alpha$ and $\theta.$  This implies that
$$
|\n |f_{\alpha \bar \beta}|^2|^2 = 0
$$
at $p$.  Since $g$ is a non-trivial solution to the differential 
equation, we conclude that the set of regular points must be dense.  
Hence the set of points at which $a_\alpha \neq 0$ for exactly 
one $\alpha$ is dense in $M.$  At such a point, let us assume 
$a_1 \neq 0$ and $a_\alpha = 0$ for all $\alpha \neq 1.$  This 
implies that
$$
\lambda_1 = -(m-1) \lambda_{\alpha}\tag 4.26
$$
and
$$
\p_{\alpha}f_{\theta \bar \theta} 
= 0 \qquad \text{for all} \qquad \alpha \neq 1. \tag 4.27
$$
Moreover, equality of \thetag{4.7} in \cite{L-W4} asserts that for
any $\theta$,
$$
\p_{\theta} f_{\alpha \bar \beta} 
= 0  \qquad \text{for all} \qquad \alpha \neq \beta 
\text{ and } \theta \neq \beta.\tag 4.28
$$
On the other hand, since
$$
\split
\p_\beta f_{\alpha \bar \beta} &= \p_\alpha f_{\beta \bar \beta}\\
&= 0 \qquad \text{if} \qquad \alpha \neq 1,
\endsplit
$$
together with \thetag {4.28}, we conclude that
$$
\p_{\theta} f_{\alpha \bar \beta} =  0  \qquad 
\text{for all} \qquad \alpha \neq \beta \text{ and } \alpha \neq 1,
\tag 4.29
$$
and
$$
\split
\p_1 f_{\theta \bar \beta} &=
\p_{\theta} f_{1 \bar \beta} \\
&= 0  \qquad \text{for all} \qquad  \beta \neq 1 
\text{ and } \theta \neq \beta.
\endsplit
\tag 4.30
$$
Also, taking the complex conjugate of \thetag {4.29} implies that
$$
\split
\bar{ \p_{\bar \theta} f_{\alpha \bar \beta}} 
&= \p_{\theta} \bar{f_{\alpha \bar \beta}}\\
&= \p_{\theta} f_{\beta \bar \alpha}\\
&=0 \qquad \text{for all} \qquad \alpha \neq \beta 
\text{ and }\alpha \neq 1.
\endsplit \tag 4.31
$$
We conclude that
$$
\n f_{\alpha \bar \beta} = 0 \qquad 
\text{for all} \qquad \alpha, \beta \neq 1 
\text{ and } \alpha \neq \beta.
$$

Moreover, the inequality for the curvature asserts that
$$
R_{\bar \alpha \beta \bar \beta \alpha} 
(\lambda_\alpha - \lambda_{\beta})^2 
\ge -\frac{\rho}{m^2}(\lambda_\alpha - \lambda_\beta)^2
$$
with equality implying that
$$
R_{\bar 1 \beta \bar\beta 1} 
=- \frac{\rho}{m^2} \qquad \text{for all} \qquad \beta \neq 1 \tag4.32
$$
since $\lambda_1 - \lambda_\beta \neq 0$ when $\beta \neq 1.$

Note that by continuity, the validity of \thetag {4.26} on a dense 
set of $M$ implies that it is valid on all of $M.$  In particular, 
by the non-vanishing of $g$, if $e_1$ is the (1,0)-vector that is 
the eigenvector for $(f_{\alpha \bar \beta})$ corresponding to the 
eigenvalue $\lambda_1$, then $e_1$ is globally defined up to 
multiplication by a complex number.  The subspace $S$ spanned 
by the other (1,0)-eigenvectors $\{e_2, \dots,e_m\}$ is also 
globally defined over $\Bbb C$.  In particular, \thetag {4.27} 
asserts that the function
$$
\split
g^2 &= |f_{\alpha \bar \beta}|\\
& = \sqrt{\frac{m}{m-1}}\,|\lambda_1|
\endsplit
$$
is constant in the directions given by $S$, hence $S$ is tangent to 
the level set of $g.$  Similarly, $\bar S$ spanned by 
$\{e_{\bar 2}, \dots, e_{\bar m}\}$ is also tangent to the level 
set of $g.$  In particular, if we define $v = \frac{\n g}{|\n g|}$, 
then the $(1,0)$-vector given by $\frac12(v -\sqrt{-1} Jv)$ must be 
a unitary multiple of $e_1.$

In any case, for $m \ge 2$, we conclude that
$$
(f_{\alpha \bar \beta}) = \left( \matrix \mu & 0 & 0 & \dots & 0\\
0 & -\frac{\mu}{m-1} & 0& \dots & 0\\
0 & 0 & -\frac{\mu}{m-1} & \dots & 0\\
\vdots & \vdots& \vdots &\,  &\vdots\\
0& 0& 0& \dots & -\frac{\mu}{m-1} \endmatrix \right) \tag 4.33
$$
where $|\mu|$ is a positive function  with
$|f_{\alpha \bar \beta}|^{\frac 12}
= \left( \frac{m}{m-1} \right)^{\frac 14}\,|\mu|^{\frac 12}$
satisfying the differential equation \thetag{4.22}.

Note that the same estimate will hold for any (real) linear
combination of such positive harmonic functions.  Indeed, this is
clear for the positive linear combination as the resulting
harmonic function is still positive. Hence, the preceding argument
works without any change. In the general case, for the sake of the
simplicity of notations, we may assume the harmonic function
$f=u-v,$ where $u$ and $v$ are two positive harmonic functions
constructed from parabolic ends $E_1$ and $E_2$ respectively. We
need to argue that the right hand side of \thetag{4.2} tends to
$0$ by choosing the cut-off function $\phi.$ In the following, we
let $\Cal L (a, b)$ denote the set $\{x\in M\,|\, a\le (u+v)(x)\le
b\}.$ On the set $E_1\cup E_2,$ we define

$$
\phi = \left\{ \aligned 1 \qquad \qquad & \qquad \text{on} \qquad
\Cal L (0, 2T) \cap (E_1\cup E_2)\\
T^{-1}(3T-u-v)& \qquad \text{on} \qquad 
\Cal L (2T, 3T) \cap (E_1\cup E_2)\\
0 \qquad \qquad & \qquad \text{on} \qquad \Cal L(3T, \infty) \cap
(E_1\cup E_2).\endaligned \right.
$$
On the remaining set $M \- (E_1\cup E_2)$, we choose $\phi= \psi\,
\chi,$ where
$$
\chi(x) = \left\{ \aligned 0 \quad \quad & \qquad 
\text{on} \qquad \Cal L(0, 2\epsilon) \\
\epsilon^{-1} (u+v - 2\epsilon)  & \qquad \text{on} 
\qquad \Cal L(2\epsilon, 3\epsilon) \\
1 \quad  \quad &  \qquad \text{on} \qquad \Cal L(3\epsilon,
\infty) \cap (M \- (E_1\cup E_2)),
\endaligned \right.
$$
and
$$
\psi(x) = \left\{ \aligned 1\quad & \qquad 
\text{on} \qquad B_\rho(R-1)\\
R-r_\rho  & \qquad \text{on} \qquad B_\rho(R) \- B_\rho(R-1)\\
0 \quad & \qquad \text{on} \qquad M\-B_\rho(R) \endaligned \right.
$$
Now the preceding argument with slight modification again shows
the right hand side of \thetag{4.2} goes to $0.$ Hence, the
equality \thetag{4.22} also holds for $f.$

Let us consider the case when $m \ge 3.$  For a fixed point 
$p \in M,$ since the complex hessian of $f \in \Cal H$ is of the 
form given by \thetag{2.33}, where $\mu_f \neq 0$ is the unique 
eigenvalue with largest absolute value, we define the sets
$$
\Cal H^+ = \{ f \in \Cal H\,|\, \mu_f >0\}
$$
and
$$
\Cal H^- = \{ f \in \Cal H \, |\, \mu_f <0\}.
$$
Obviously, both $\Cal H^+$ and $\Cal H^-$ are nonempty because 
$f\in \Cal H^+$ implies that $-f \in \Cal H^-.$  Moreover, 
$\Cal H^+ \cap \Cal H^- = \emptyset$ and 
$\Cal H \-\{0\} = \Cal H^+ \cup \Cal H^-.$  This is only possible 
if $\dim \Cal H = 1,$ hence $M$ has at most 2 ends. 

When $m=2$, the above argument is not valid since $\mu_f$ is not 
unique.  However, we will show that $M$ must have at most 4 ends 
in this case. Indeed, for a fixed point $p \in M$, let us define 
the map sending $f \in \Cal H$ to
the unitary eigenvector $e_1$ of $(f_{\alpha \bar \beta})$ 
corresponding to the unique positive eigenvalue $\mu >0.$  
Since $e_1$ is defined only up to complex scalar multiplication 
and the map is invariant under positive scalar multiplication of $f$, 
it induces a map
$$
F: \Bbb S(\Cal H) \longrightarrow \Bbb{CP}^1(T_pM)
$$
where $\Bbb S (\Cal H)$ is the unit sphere of the vector space 
$\Cal H$ and $\Bbb{CP}^1(T_pM)$ is the set of complex lines in the 
tangent space of $M$ at $p.$  If  $\dim \Cal H \ge 4,$ 
then $\dim \Bbb S(\Cal H) \ge 3.$  A simple dimension argument 
implies that the map $F$ cannot be injective.  Hence there exists 
two harmonic functions $f_1$ and $f_2$ with the same eigenvector 
with eigenvalues $\mu_1$ and $\mu_2$.  After taking an appropriate 
linear combination of $f_1$ and $f_2$, we produce a non-trivial 
harmonic function in $\Cal H$ with vanishing complex Hessian at $p.$  
This gives a contradiction, hence $M$ cannot have more than $4$ ends.

When $M$ has exactly $4$ ends and $\dim \Cal H = 3$, the  
map $F$ maps a $2$-dimensional sphere $\Bbb S(\Cal H)$ injectively 
into $\Bbb {CP}^1(T_pM).$  Since $\Bbb {CP}^1$ is homeomorphic 
to $\Bbb S^2$, this implies that $F$ must be onto.  In particular, 
any (1,0)-vector $e_1$ in $T_pM$ can be realized as an eigenvector 
of the complex Hessian of some $f\in \Cal H.$  Using \thetag{4.32}, 
we conclude that all holomorphic bisectional curvatures must be 
given by $-\rho.$  Since the point $p$ is arbitrary, this shows that
$$
R_{\bar \alpha \beta \bar \beta \alpha}(x) = -\frac{\rho(x)}4
$$
for all $\alpha \neq \beta.$ \qed
\enddemo

\enddocument